\documentclass[12pt, a4paper, reqno]{amsart}
\usepackage{preemble}
\usepackage{microtype}
\begin{document}

\title[Dimension conservation of harmonic measures]
      {Dimension conservation of harmonic measures in products of hyperbolic spaces}
      
      \author{Ryokichi Tanaka}
      \address{Department of Mathematics, 
      Kyoto University, Kyoto 606-8502 JAPAN}
      \email{rtanaka@math.kyoto-u.ac.jp}
\date{\today}

\maketitle

\begin{abstract}
We show that the harmonic measure on a product of boundaries satisfies dimension conservation for a random walk with non-elementary marginals on a countable group acting on a product of hyperbolic spaces under the finite first moment condition.
\end{abstract}

\section{Introduction}

Let $\Gamma$ and $\Gammast$ be non-elementary hyperbolic groups.
We study a random walk on the product group $\Gammab:=\Gamma \times \Gammast$ and establish a dimension formula for the harmonic measure on the product of (Gromov) boundaries.
After stating our results in this special case, we consider a countable group of isometries of a product of two hyperbolic metric spaces.

Let $\pi$ be a probability measure on $\Gammab$ such that marginals $\mu$ and $\must$ are non-elementary, i.e., their supports generate non-elementary subgroups in $\Gamma$ and in $\Gammast$ as groups respectively.
For such a $\pi$, a {\bf harmonic measure} $\nu_\pi$ is defined on the product of boundaries $\partial \Gamma \times \partial \Gammast$ (cf.\ Section \ref{Sec:hyperbolic}).
It is the unique probability measure satisfying that
\[
\nu_\pi=\pi\ast \nu_\pi, \quad \text{where $\pi\ast \nu_\pi=\sum_{\xb \in \Gammab}\pi(\xb)\xb\nu_\pi$ and $\xb\nu_\pi:=\nu_\pi\circ \xb^{-1}$}.
\]
In the above, we consider the natural action of $\Gammab=\Gamma\times \Gammast$ on $\partial \Gamma\times \partial \Gammast$.
The harmonic measure $\nu_\pi$ has marginals $\nu_\mu$ and $\nu_{\must}$ on $\partial \Gamma$ and on $\partial \Gammast$ respectively, and these are determined by $\mu$ and by $\must$.
The measure $\nu_\mu$ (or $\nu_{\must}$) for a single hyperbolic group and its generalization has attracted intensive studies, including the dimension, e.g., \cite{KaimanovichTrees, LedrappierFreeGroups}, more recently, \cite{BHM11, HochmanSolomyak, Tdim, DussauleYangHausdorff}.
However, the harmonic measure $\nu_\pi$ for a product group has been studied only in a few cases (see \cite{VolkovThesis} for a special case of products of hyperbolic free product groups).
Dimensional properties of such a measure exhibit new features since it contains different behaviors depending on the factors.
This manifests an additional difficulty, which arises in a higher rank setting-----in that case, the boundaries are assembled in further intricate ways such as flag varieties.
See \cite{KaimanovichLePrince} for a thorough discussion on this matter of subject, and recent works \cite{LessaDisintegrations, RapaportExact, LedrappierLessaExact}.
For related results on ergodic invariant measures for affine iterated function systems,
see e.g., \cite{Kenyon-Peres, Feng-Hu, Feng-dimension}.
We study the harmonic measure $\nu_\pi$ itself and the conditional measure $\nu_\pi^\eta$ of $\nu_\pi$ on $\partial \Gamma\times \partial \Gammast$ for $\nu_{\must}$-almost every $\eta \in \partial \Gammast$.
It is shown that these are exact dimensional under a finite first moment condition.
The dimension formula is a sum of ratios of asymptotic entropy over drift corrected with differential entropy.
For the conditional measures, we provide a sufficient condition for the strict positivity of dimension.

Let us consider left invariant hyperbolic metrics $d$ and $\dst$ quasi-isometric to word metrics in $\Gamma$ and in $\Gammast$, respectively.
They induce quasi-metrics $q$ and $\qst$ in the compactified spaces $\Gamma\cup\partial \Gamma$ and $\Gammast\cup\partial \Gammast$, respectively, and let $\qb(\xib_1, \xib_2):=\max\{q(\xi_1, \xi_2), \qst(\eta_1, \eta_2)\}$ for $\xib_i=(\xi_i, \eta_i) \in (\Gamma\cup\partial \Gamma)\times (\Gammast\cup\partial \Gammast)$ and $i=1, 2$ (cf.\ Section \ref{Sec:hyperbolic}).
Let us assume that $\pi$ has a finite first moment, i.e.,
$\sum_{\xb \in \Gammab}\db(\idb, \xb)\pi(\xb)<\infty$ for the identity element $\idb$ and $\db(\xb_1, \xb_2):=\max\{d(x_1, x_2), \dst(y_1, y_2)\}$ for $\xb_i=(x_i, y_i)\in \Gammab$ and $i=1, 2$.
The {\bf asymptotic entropy} $h(\pi)$ is defined as the limit 
for $n$-fold convolutions $\pi_n:=\pi^{\ast n}$,
\[
h(\pi)=\lim_{n \to \infty}\frac{1}{n}\sum_{\xb \in \Gammab}-\pi_n(\xb)\log \pi_n(\xb).
\]
The {\bf drift} $l(\Gamma, \mu)$ is defined as the limit for $\mu_n:=\mu^{\ast n}$,
\[
l(\Gamma, \mu)=\lim_{n \to \infty}\frac{1}{n}\sum_{x \in \Gamma}d(\id, x)\mu_n(x).
\]
Similarly for $l(\Gammast, \must)$.
These are positive since $\mu$ and $\must$ are non-elementary (cf.\ Section \ref{Sec:RW}).
For every real $r > 0$ and 
and every $\xib \in \partial \Gamma\times \partial \Gammast$,
let $\Bb(\xib, r)$ denote the open ball of radius $r$ centered at $\xib$ in $(\partial \Gamma\times \partial \Gammast, \qb)$.

\begin{theorem}\label{Thm:product_intro}
Let $\Gamma$ and $\Gammast$ be non-elementary hyperbolic groups, and $\pi$ be a probability measure on $\Gammab=\Gamma\times \Gammast$ with finite first moment and non-elementary marginals $\mu$ and $\must$ on $\Gamma$ and on $\Gammast$ respectively.
Suppose that $l(\Gamma, \mu) \ge l(\Gammast, \must)$.
Then
the harmonic measure $\nu_\pi$ on $\(\partial \Gamma \times \partial \Gammast, \qb\)$ is exact dimensional,
i.e., for $\nu_\pi$-almost every $\xib \in \partial \Gamma \times \partial \Gammast$,
\[
\lim_{r \to 0}\frac{\log \nu_\pi\(\Bb(\xib, r)\)}{\log r}=\frac{h(\pi)-h(\must)}{l(\Gamma, \mu)}+\frac{h(\must)}{l(\Gammast, \must)}.
\]
In particular, the Hausdorff dimension of $\nu_\pi$ is computed as
\[
\dim \nu_\pi= \frac{h(\pi)-h(\must)}{l(\Gamma, \mu)}+\frac{h(\must)}{l(\Gammast, \must)}.
\]
\end{theorem}

Theorem \ref{Thm:product_intro} is shown in Theorem \ref{Thm:product} in a more general setting.
The above result is based on the exact dimensionality of disintegrated measures:
Let $\nu_\pi^\eta$ for $\eta \in \partial \Xcst$ denote a system of conditional measures of $\nu_\pi$ on $\partial \Gamma \times \partial \Gammast$ with respect to the $\sigma$-algebra generated by the projection from $\partial \Gamma\times \partial \Gammast$ to $\partial \Gammast$.

\begin{theorem}\label{Thm:exact_intro}
Let $\Gamma$ and $\Gammast$ be non-elementary hyperbolic groups, and $\pi$ be a probability measure on $\Gammab=\Gamma\times \Gammast$ with finite first moment and non-elementary marginals $\mu$ and $\must$ on $\Gamma$ and on $\Gammast$ respectively.
For $\nu_{\must}$-almost every $\eta \in \partial \Gammast$,
the conditional measure $\nu_\pi^\eta$ is exact dimensional on $(\partial \Gamma\times \partial \Gammast, \qb)$,
i.e., for $\nu_\pi^\eta$-almost every $\xib \in \partial \Gamma\times \partial \Gammast$,
\[
\lim_{r \to 0}\frac{\log \nu_\pi^\eta(\Bb(\xib, r))}{\log r}=\frac{h(\pi)-h(\must)}{l(\Gamma, \mu)},
\]
In particular, 
the Hausdorff dimension of $\nu_\pi^\eta$ is computed as 
for $\nu_{\must}$-almost every $\eta \in \partial \Gammast$,
\[
\dim \nu_\pi^\eta=\frac{h(\pi)-h(\must)}{l(\Gamma, \mu)}.
\]
\end{theorem}

Theorem \ref{Thm:exact_intro} is shown in Theorem \ref{Thm:exact} in a more general setting.
Following Furstenberg \cite[Definition 3.1]{Furstenberg-Ergodic},
we say that a Borel probability measure $\nu$ on a product of compact metric spaces $\Mcc \times \Mccst$
satisfies {\bf dimension conservation} if the following holds.
Let us consider the pushforward $\nust$ and a system of conditional measures $\nu^\eta$ for $\eta \in \Mccst$ of $\nu$ associated with the projection $\Mcc \times \Mccst \to \Mccst$: 
the measures
$\nu$ and $\nust$ are exact dimensional with dimension $\dim \nu$ and $\dim \nust$ respectively,
for $\nust$-almost every $\eta \in \Mccst$, conditional measures $\nu^\eta$ are exact dimensional with dimension $\dim \nu^\eta$, and
\[
\dim \nu=\dim \nu^\eta+\dim \nust.
\]
In the above, we understand that the metric in $\Mcc \times \Mccst$ is defined as the maximum of metrics in $\Mcc$ and $\Mccst$ (or an arbitrary one bi-Lipschitz to it).
It has been shown that $\nu_\must$ on $(\partial \Gammast, \qst)$ is exact dimensional with dimension $h(\must)/l(\Gammast, \must)$ \cite[Theorem 3.8]{Tdim}.
Therefore
Theorems \ref{Thm:product_intro} and \ref{Thm:exact_intro} imply that
the harmonic measure $\nu_\pi$ on $\partial \Gamma \times \partial \Gammast$ satisfies dimension conservation.
The statement holds in a more general setting; see Section \ref{Sec:product}.
For an extension to a product of more than two hyperbolic groups, see Remark \ref{Rem:two_or_more}.

Let $(\Xc, d)$ and $(\Xcst, \dst)$ be proper roughly geodesic hyperbolic metric spaces with bounded growth at some scales (for the definitions, see Section \ref{Sec:pre}).
Examples of such spaces include Gromov hyperbolic Riemannian manifolds with sectional curvature bounded from below and from above and Cayley graphs of hyperbolic groups.
The space $\Xc \times \Xcst$ is equipped with a base point $\ob$ and the metric $\db(\xb_1, \xb_2):=\max\{d(x_1, x_2), \dst(y_1, y_2)\}$ for $\xb_i=(x_i, y_i) \in \Xc\times \Xcst$ and $i=1, 2$.
Let us consider a countable subgroup $\Gammab$ in the product of isometry groups $\Isom \Xc \times \Isom \Xcst$.
We say that $\Gammab$ has a {\bf finite exponential growth} relative to $(\Xc \times \Xcst, \db)$
if there exists a constant $c>0$ such that for all $r > 0$,
\[
\#\big\{\xb \in \Gammab \ : \ \db(\ob, \xb \cdot \ob)< r\big\}\le ce^{c r}.
\]
In the above, $\#A$ denotes the cardinality of a set $A$.
For a probability measure $\pi$ on $\Gammab$, let $\mu$ and $\must$ denote the marginal on $\Isom \Xc$ and on $\Isom \Xcst$ respectively.
Let $\supp \must$ denote the support of $\must$.
The {\bf differential entropy} of the pair $(\partial \Xcst, \must)$ is defined by
\[
h(\partial \Xcst, \must):=\sum_{x\in \supp \must}\must(x)\int_{\partial \Xcst}\log \frac{dx\nu_{\must}}{d\nu_{\must}}(\eta)\,dx\nu_{\must}(\eta).
\]
In general, it holds that $h(\partial \Xcst, \must) \le h(\must)$,
and the equality holds if and only if $(\partial \Xcst, \nu_\must)$ is a Poisson boundary for the pair $(\Isom \Xcst, \must)$ (cf.\ Section \ref{Sec:RW}).
Let $l(\Xc, \mu)$ be the drift associated with a $\mu$-random walk on $\Xc$.
Theorem \ref{Thm:exact_intro} is generalized in this setting.

\begin{theorem}\label{Thm:exact}
Let $(\Xc, d)$ and $(\Xcst, \dst)$ be proper roughly geodesic hyperbolic metric spaces with bounded growth at some scale, and $\Gammab$ be a countable subgroup of $\Isom \Xc\times \Isom \Xcst$ with finite exponential growth relative to $(\Xc\times \Xcst, \db)$.
If $\pi$ is a probability measure on $\Gammab$ with finite first moment and non-elementary marginals $\mu$ and $\must$ respectively,
then the conditional measure $\nu_\pi^\eta$ is exact dimensional for $\nu_{\must}$-almost every $\eta \in \partial \Xcst$.
In fact, for $\nu_{\must}$-almost every $\eta \in \partial \Xcst$ and for $\nu_\pi^\eta$-almost every $\xib \in \partial \Xc\times \partial \Xcst$,
\[
\lim_{r \to 0}\frac{\log \nu_\pi^\eta(\Bb(\xib, r))}{\log r}=\frac{h(\pi)-h(\partial \Xcst, \must)}{l(\Xc, \mu)}.
\]
In particular, the Hausdorff dimension of $\nu_\pi^\eta$ is computed as for $\nu_{\must}$-almost every $\eta \in \partial \Xcst$,
\[
\dim \nu_\pi^\eta=\frac{h(\pi)-h(\partial \Xcst, \must)}{l(\Xc, \mu)}.
\]
\end{theorem}

For the differential entropy, 
it holds that $h(\partial \Xcst, \must)=0$ if and only if $(\partial \Xcst, \nu_{\must})$ is trivial, i.e., $\nu_{\must}$ is invariant under the action of $\Gammab$ on $\partial \Xcst$.
If $\must$ is non-elementary,
then $(\partial \Xcst, \nu_{\must})$ is non-trivial and $h(\partial \Xcst, \must)>0$.
In the setting of Theorem \ref{Thm:exact},
it can be the case that $h(\pi)=h(\partial \Xcst, \must)$ (see Example \ref{Ex:ns} below).
The following result provides a sufficient condition under which $h(\pi)>h(\partial \Xcst, \must)$, i.e., $(\partial \Xcst, \nu_\must)$ is a proper quotient of the Poisson boundary for the pair $(\Gammab, \pi)$.
If this is the case, then the Hausdorff dimension of conditional measures are strictly positive.

\begin{theorem}\label{Thm:entropy_lower_bound_intro}
Let $\Gamma$ and $\Gammast$ be countable subgroups in $\Isom \Xc$ and in $\Isom \Xcst$ respectively, and $\Gammab:=\Gamma\times \Gammast$.
Further let us consider a probability measure $\pi$ on $\Gammab$ of the following form:
For some $\alpha \in (0, 1]$, 
\[
\pi=\alpha \lambda \times \lambdast+(1-\alpha)\pi_0
\]
with non-elementary probability measures $\lambda$ and $\lambdast$ on $\Gamma$ and on $\Gammast$ respectively, and a probability measure $\pi_0$ on $\Gammab$.
It holds that $h(\pi)-h(\partial \Xcst, \must)>0$,
where $\must$ is the marginal of $\pi$ on $\Gammast$.
\end{theorem}

Theorem \ref{Thm:entropy_lower_bound_intro} is shown in Theorem \ref{Thm:entropy_lower_bound};
moreover, if in addition $\Gammab=\Gamma \times \Gammast$ has a finite exponential growth relative to $(\Xc\times \Xcst, \db)$,
then
for $\nu_{\must}$-almost every $\eta \in \partial \Xcst$,
the Hausdorff dimension of the conditional measure $\nu_\pi^\eta$ is positive (cf.\ Theorem \ref{Thm:exact}).

\begin{example}\label{Ex:ns}
Let $\Gamma$ be a hyperbolic group and $\mu$ be a non-elementary probability measure on $\Gamma$ with finite first moment relative to a word metric.
For $\rho \in [0, 1]$,
let
\[
\pi^\rho:=\rho \mu\times \mu+(1-\rho)\mu_{\diag},
\]
where $\mu_{\diag}((x, \xst)):=\mu(x)$ if $x=\xst$, and $0$ if otherwise.
The $\pi^\rho$-random walk on $\Gamma \times \Gamma$ appears in the study of noise sensitivity problem on groups \cite{BenjaminiBrieussel, NonNS}.
By Theorems \ref{Thm:exact_intro} and \ref{Thm:entropy_lower_bound_intro} applied to the case when $\Gamma=\Gammast$ and $\mu=\must$,
it holds that for all $\rho \in (0, 1]$,
\[
\dim \nu_{\pi^\rho}^\eta=\frac{h(\pi^\rho)-h(\mu)}{l(\Gamma, \mu)}>0 \quad \text{for $\nu_\mu$-almost every $\eta \in \partial \Gamma$}.
\]
For $\rho=0$, since $h(\pi^\rho)=h(\mu)$, it holds that $\dim \nu_{\pi^\rho}^\eta=0$ for $\nu_\mu$-almost every $\eta \in \partial \Gamma$.
Theorem \ref{Thm:product_intro} shows that for all $\rho \in [0, 1]$,
\[
\dim \nu_{\pi^\rho}=\frac{h(\pi^\rho)}{l(\Gamma, \mu)}.
\]
This reproduces \cite[Theorem 3.1]{NonNS}.
\end{example}

\begin{example}\label{Ex:diagonal}
This example is not covered by Theorem \ref{Thm:entropy_lower_bound_intro}.
Suppose that $\Gammab=\Gamma \times \Gammast$ for two hyperbolic groups $\Gamma$ and $\Gammast$, and
that there exists a (non-injective) surjective homomorphism $\Pi:\Gamma \to \Gammast$.
Let $\Delta:\Gamma\to \Gamma\times \Gammast$ be the diagonal embedding $\Delta(x)=(x, \Pi(x))$ for $x \in \Gamma$.
For a non-elementary probability measure $\mu$ on $\Gamma$ with finite first moment,
let $\pi:=\Delta_\ast \mu$ be the pushforward of $\mu$ by $\Delta$.
In this case, marginals of $\pi$ are $\mu$ on $\Gamma$ and $\Pi_\ast \mu$ on $\Gammast$ respectively,
where $\Pi_\ast \mu$ is the pushforward of $\mu$ by $\Pi$.
Applying to Theorem \ref{Thm:exact_intro} with $\must=\Pi_\ast \mu$ shows that
\[
\dim \nu_\pi^\eta=\frac{h(\mu)-h(\Pi_\ast \mu)}{l(\Gamma, \mu)} \quad \text{for $\nu_{\Pi_\ast \mu}$-almost every $\eta \in \partial \Gammast$}.
\]
This follows since $h(\pi)=h(\mu)$ and $h(\partial \Gammast, \Pi_\ast \mu)=h(\Pi_\ast \mu)$.
It holds that $h(\mu)=h(\Pi_\ast \mu)$ if $\Pi$ is an isomorphism, and
it depends on $\Pi$ whether a strict inequality $h(\mu)>h(\Pi_\ast \mu)$ holds or not.
As a simple explicit example, 
let $\Gamma=F_{m+1}$ and $\Gammast=F_m$ be free groups of rank $m+1$ and $m$ respectively for $m \ge 2$, equipped with word metrics associated with free bases.
Further let $\Pi:F_{m+1} \to F_m$ be a homomorphism sending $x_i$ to $y_i$, $i=1, \dots, m$, and sending $x_{m+1}$ to the identity, where $\{x_1, \dots, x_{m+1}\}$ denotes the free basis in $F_{m+1}$ and $\{y_1, \dots, y_m\}$ denotes the free basis in $F_m$.
For the uniform distribution $\mu$ on the symmetrized free basis in $F_{m+1}$,
the induced distribution on $F_{m}$ defines a simple random walk on $F_m$ with holding probability $1/(m+1)$.
A computation yields $l(F_{m+1}, \mu)=m/(m+1)$, $l(F_m, \Pi_\ast \mu)=(m-1)/(m+1)$,
\[
h(\mu)=\frac{m}{m+1}\log(2m+1) \quad \text{and }\quad h(\Pi_\ast \mu)=\frac{m-1}{m+1}\log(2m-1).
\]
Therefore 
\[
\dim \nu_\pi^\eta=\log(2m+1)-\frac{m-1}{m}\log (2m-1) \quad \text{for $\nu_\pi^\eta$-almost every $\eta \in \partial F_m$}.
\]
Furthermore, Theorem \ref{Thm:product_intro} shows that
\[
\dim \nu_\pi=\log(2m+1)+\frac{1}{m}\log(2m-1).
\]
\end{example}

\subsection*{Outlines of proofs}
Let us briefly mention the proof of Theorem \ref{Thm:product_intro} for a product of two hyperbolic groups.
For a single hyperbolic group $\Gamma$ with a non-elementary probability measure $\mu$,
the corresponding harmonic measure $\nu_\mu$ on $\partial \Gamma$ is exact dimensional \cite[Theorem 3.8]{Tdim}.
Roughly speaking, it boils down to estimate probabilities that for a $\mu$-random walk $w_n$, an independent $\mu$-random walk $w_n'$ is around $w_n$ within distance $o(n)$ for $n \in \Z_+$.
This leads an estimate of the harmonic measure $\nu_\mu$ on the balls $B(w_\infty, e^{-ln})$ where $l:=l(\Gamma, \mu)$. 
Here the $\mu$-random walk $\{w_n\}_{n \in \Z_+}$ is for a sampling $w_\infty$ in $\partial \Gamma$ according to $\nu_\mu$ and the independent $\mu$-random walk $\{w_n'\}_{n \in \Z_+}$ is for the estimate $\nu_\mu(B(w_\infty, e^{-ln}))$.
Since the probability that $w_n'$ is around $w_n$ within distance $o(n)$ is $e^{-h(\mu)n+o(n)}$
by the Shannon theorem for random walks, 
this explains $\nu_\mu(B(w_\infty, e^{-ln}))=e^{-h(\mu)n+o(n)}$, which is the exact dimensionality of $\nu_\mu$ with the right dimension $h(\mu)/l$.

The conditional measure $\nu_\pi^\eta$ for $\nu_\must$-almost every $\eta \in \partial \Gammast$ is the hitting distribution of a conditional process. 
This is a Markov chain (although the transition probabilities are not group-invariant) and (one of) the methods developed for a single hyperbolic group in [ibid] applies.
The asymptotic entropy of this conditional process equals $h(\pi)-h(\must)$ by the Shannon theorem for the conditional process \cite{Kaimanovich-hyperbolic}.
The conditional measure $\nu_\pi^\eta$ is defined on $\partial \Gamma \times \partial \Gammast$ but supported on $\partial \Gamma \times \{\eta\}$ for $\nu_\must$-almost every $\eta \in \partial \Gammast$.
An analogous discussion to the $\mu$-random walk above works and this leads to estimating the $\nu_\pi^\eta$-measures on the balls in the boundary $B(w_\infty, e^{-l n})\times \{\eta\}$ for $\nu_\pi^\eta$-almost every $(w_\infty, \eta) \in \partial \Gamma \times \partial \Gammast$.
In fact, we obtain $\nu_\pi^\eta(B(w_\infty, e^{-ln})\times \{\eta\})=e^{-(h(\pi)-h(\must))n+o(n)}$, deducing Theorem \ref{Thm:exact_intro}.

The harmonic measure $\nu_\pi$ on $\partial \Gamma \times \partial \Gammast$ is, however, analyzed in a completely different way.
First of all it requires to take into account the difference between $l$ and $\lst$ where $\lst:=l(\Gammast, \must)$.
If $l \ge \lst$, then
\[
\frac{h(\pi)-h(\mu)}{\lst}+\frac{h(\mu)}{l}\le \frac{h(\pi)-h(\must)}{l}+\frac{h(\must)}{\lst},
\]
since $h(\pi) \le h(\mu)+h(\must)$, and the inequality can be strict.
Since the right hand side of the above inequality is the correct value, the dimension upper bound should use the inequality $l\ge \lst$ whereas the dimension lower bound would not need it.
Concerning the dimension upper bound, the Shannon theorem for the conditional process shows that for $\nu_\must$-almost every $\eta \in \partial \Xcst$, for the conditional process $\{\wb_n\}_{n \in \Z_+}$,
\[
\Pb^\eta([\wb_0, \dots, \wb_n]) \le e^{h(\must)n+o(n)}\Pb([\wb_0, \dots, \wb_n]).
\]
In the above, $[\wb_0, \dots, \wb_n]$ denotes the cylinder set.
At this point, we keep track the whole trajectory up to time $n$ instead of just looking at the position $\wb_n$.
The argument here is inspired by \cite[Section 8]{LedrappierLessaExact} (where they refer to \cite{Feng-dimension} for the idea).
The $\nu_\pi^\eta$ on the balls $\Bb(\wb_\infty, e^{-\lst n})=B(w_\infty, e^{-\lst n})\times B(\wst_\infty, e^{-\lst n})$ estimates
by Theorem \ref{Thm:exact_intro},
\[
\nu_\pi^\eta\(\Bb(\wb_\infty, e^{-\lst n})\)=\exp\(-\(\frac{h(\pi)-h(\must)}{l}\)\lst n+o(n)\).
\]
Averaging $\eta$ over $B(\wst_\infty, e^{-\lst n})$ deduces the required lower bound (thus upper bound for the dimension) of $\nu_\pi(\Bb(\wb_\infty, e^{-\lst n}))$.
In this discussion, it is crucial to use the balls with radii $e^{-\lst n}$ rather than $e^{-l n}$ (or other scales)
since $\qb(\wb_\infty, \wb_n)=e^{-\lst n+o(n)}$, where $l \ge \lst$,
\[
q(w_\infty, w_n)=e^{-ln+o(n)} \quad \text{and} \quad \qst(\wst_\infty, \wst_n)=e^{-\lst n+o(n)}.
\]
Concerning the dimension lower bound, a slight strengthened version for the lower bound in Theorem \ref{Thm:exact_intro} enables us to exploit the naive disintegration formula.
Roughly, estimating along the following heuristic can be justified:
\[
\nu_\pi\(\Bb(\wb_\infty, e^{- ln})\)
\approx \nu_\pi^\eta (B(w_\infty, e^{-ln})\times \{\eta\}) \cdot \nu_\must\(B(\wst_\infty, e^{-ln})\).
\]
Since this works only for $\nu_\pi$ restricted on a large subset, the argument is merely for the upper bound (thus lower bound for the dimension) of $\nu_\pi$ (up to a density lemma which is guaranteed by a weak version of the Lebesgue differentiation theorem Lemma \ref{Lem:BonkSchramm}).
Thus Theorem \ref{Thm:exact_intro} and the exact dimensionality of $\nu_\must$ with $h(\must)/\lst$ conclude the required dimension lower bound on $\nu_\pi$.

The above sketch for hyperbolic groups can be extended to a countable group of isometries acting on a product of two hyperbolic metric spaces in Theorem \ref{Thm:product}.
The positive lower bound for $h(\pi)-h(\partial \Xcst, \must)$ in Theorem \ref{Thm:entropy_lower_bound_intro} uses the pivotal time technique developed by Gou\"ezel in \cite{Gouezel_Exp}.
We mention possible extensions of Theorems \ref{Thm:product_intro} and \ref{Thm:product} and questions in Remark \ref{Rem:two_or_more}.

\subsection*{Organization}
Section \ref{Sec:pre} recalls basics on hyperbolic metric spaces and random walks.
Section \ref{Sec:dimension} concerns dimensions of the conditional measures, showing Theorem \ref{Thm:exact} (and thus Theorem \ref{Thm:exact_intro}).
Section \ref{Sec:product} concerns dimensions of the harmonic measures on products of boundaries, showing Theorem \ref{Thm:product} (and thus Theorem \ref{Thm:product_intro}). 
Section \ref{Sec:positive} is about a sufficient condition on a positivity of the dimension for conditional measures, showing Theorem \ref{Thm:entropy_lower_bound_intro} in Theorem \ref{Thm:entropy_lower_bound}. 

\subsection*{Notation}
We denote by $c$, $C$, \dots, constants whose exact values may vary from line to line, and by $C_\delta$ a constant which depends on the other constant $\delta$ to emphasize its dependency.
For a real valued sequence $\{f(n)\}_{n \in \Z_+}$ on non-negative integers $\Z_+$, we write $f(n)=o(n)$ if $|f(n)|/n \to 0$ as $n \to \infty$.
For a set $A$, we denote by $A^{\sf c}$ the complement set, and by $\#A$ the cardinality.

\section{Preliminaries}\label{Sec:pre}

\subsection{Hyperbolic metric spaces}\label{Sec:hyperbolic}

For background, we refer to the original paper by Gromov \cite{GromovHyperbolic}.
For a metric space $(\Xc, d)$, the {\bf Gromov product} is defined by
\[
(x|y)_z:=\frac{1}{2}(d(x, z)+d(z, y)-d(x, y)) \quad \text{for $x, y, z \in \Xc$}.
\]
A metric space $(\Xc, d)$ is {\bf $\delta$-hyperbolic} for a non-negative real $\delta \in \R_+$
if it holds that
\begin{equation}\label{Eq:delta}
(x|y)_w \ge \min\big\{(x|z)_w, (z|y)_w\big\}-\delta \quad \text{for all $x, y, z, w \in \Xc$}.
\end{equation}
It is called {\bf hyperbolic} if it is $\delta$-hyperbolic for some $\delta \in \R_+$.
A map $\gamma: I \to \Xc$ from an interval $I$ in $\R$ to $\Xc$ is called a $C$-{\bf rough geodesic} for $C \in \R_+$ if $|d(\gamma(s), \gamma(t))-|t-s||\le C$ for all $s, t \in I$.
Further a map $\gamma: I \to \Xc$ is called a $C$-{\bf rough geodesic ray} in the case when $I=[0, \infty)$.
A metric space is called $C$-{\bf roughly geodesic} for $C \in \R_+$ if for all pairs of points $x, y \in \Xc$ there exists a $C$-rough geodesic $\gamma: [a, b] \to \Xc$ such that $\gamma(a)=x$ and $\gamma(b)=y$.
In this terminology, a metric space is called {\bf geodesic} if it is $0$-roughly geodesic.
A graph endowed with a path metric of unit edge length (e.g., a Cayley graph) is also considered as a geodesic metric space by using intervals in the integers $\Z$ in the definition.
Let us simply call a metric space {\bf roughly geodesic} if it is $C$-roughly geodesic for some $C \in \R_+$.
For a hyperbolic group $\Gamma$ equipped with a left invariant hyperbolic metric $d$ quasi-isometric to a word metric, $(\Gamma, d)$ is roughly geodesic (cf.\ \cite[Proposition 5.6]{BonkSchramm} and \cite[Theorem 2.2]{BHM11}).
A metric space $(\Xc, d)$ is {\bf proper} if for all $x \in \Xc$ and all $r\in \R_+$,
the ball
$B(x, r):=\big\{y \in \Xc \ : \ d(x, y) < r\big\}$
is relatively compact.

For a hyperbolic metric space $(\Xc, d)$,
the {\bf (Gromov) boundary} $\partial \Xc$ is defined as the set of equivalence classes of divergent sequences in $\Xc$.
Let us fix a point $o \in \Xc$.
Further let $q(x, y):=\exp(-(x|y)_o)$ for $x, y \in \Xc$ and $q(x,y):=0$ if $x=y$.
Since the space is $\delta$-hyperbolic for some $\delta \in \R_+$,
it holds that
\begin{equation}\label{Eq:q}
q(x, y)\le e^\delta \max\big\{q(x, z), q(z, y)\big\} \quad \text{for $x, y, z \in \Xc$}.
\end{equation}
A sequence $\{x_n\}_{n \in \Z_+}$ in $\Xc$ is called {\bf divergent} if it is a Cauchy sequence with respect to $q$.
Two sequences $\{x_n\}_{n \in \Z_+}$ and $\{y_n\}_{n \in \Z_+}$ are {\bf equivalent} if $q(x_n, y_m) \to 0$ as $n, m \to \infty$.
It is indeed an equivalence relation in the set of divergence sequences by \eqref{Eq:q}.

For $\xi \in \Xc \cup \partial \Xc$,
let us write $\xi=[\{x_n\}_{n \in \Z_+}]$ for a divergent sequence $\{x_n\}_{n \in \Z_+}$ which represents $\xi$ if $\xi \in \partial \Xc$, 
or for the constant sequence $x_n=\xi$ for all $n \in \Z_+$ if $\xi \in \Xc$.
The Gromov product is extended to $\Xc \cup \partial \Xc$ by 
\[
(\xi|\eta)_o:=\inf\Big\{\liminf_{n, m \to \infty}(x_n|y_m)_o \ : \ \xi=[\{x_n\}_{n \in \Z_+}], \ \eta=[\{y_m\}_{m \in \Z_+}]\Big\}.
\]
For a $\delta$-hyperbolic space, 
the extended Gromov product satisfies \eqref{Eq:delta} for $x, y, z \in \Xc \cup \partial \Xc$ and $w=o$.
Let us extend $q$ on $\Xc$ to $\Xc \cup \partial \Xc$ and call it the {\bf quasi-metric}: 
\[
q(\xi, \eta):=\exp(-(\xi|\eta)_o) \quad \text{if $\xi\neq \eta$}, \quad 
\text{and} \quad q(\xi, \eta):=0 \quad \text{if $\xi=\eta$}, \quad \text{for $\xi, \eta \in \Xc \cup \partial \Xc$}.
\]
It is known that there exists an $\e_0>0$ such that for all $0<\e<\e_0$ the power $q^\e$ is bi-Lipschitz equivalent to a genuine metric.
However, the quasi-metric $q$ is used to define balls and other notions related to metrics without introducing an additional parameter $\e$.
The space $\Xc \cup \partial \Xc$ is equipped with the topology defined from the (quasi-)metric.
If $\Xc$ is proper, then $\Xc \cup \partial \Xc$ is a compact metrizable space.
If in addition $\Xc$ is $C$-roughly geodesic for some $C \in \R_+$,
then for every $\xi \in \partial \Xc$ and every $x \in \Xc$ there exists a $C$-rough geodesic ray $\gamma$ from $x$ converging to $\xi$, i.e., $q(\gamma(t), \xi) \to 0$ as $t \to \infty$ in $\Xc \cup \partial \Xc$ (cf.\ \cite[Proposition 5.2]{BonkSchramm}).
Henceforth it is assumed that $\Xc$ is a proper roughly geodesic hyperbolic metric space.

Let us denote the open ball of radius $r\in\R_+$ centered at $\xi \in \Xc\cup\partial \Xc$ in $\Xc\cup \partial \Xc$ by 
\[
B(\xi, r):=\Big\{\eta \in \Xc\cup\partial \Xc \ : \ q(\xi, \eta)< r\Big\}.
\]
The {\bf shadow} (seen from $o$) at $x \in \Xc$ with thickness $R \in \R_+$ is defined by
\[
\Oc(x, R):=\Big\{\eta \in \partial \Xc \ : \ (o|\eta)_x < R\Big\}.
\]
The following is used to compare shadows with balls. 
For each $T>0$,
there exist constants $R_0, C>0$ such that
for all $R>R_0$, all $\xi \in \partial \Xc$ and all $x \in \Xc$ with $(o|\xi)_x \le T$,
\begin{equation}\label{Eq:shadows_balls}
B(\xi, C^{-1}e^{-d(o, x)+R})\cap \partial \Xc \subset \Oc(x, R) \subset B(\xi, C e^{-d(o, x)+R})\cap \partial \Xc.
\end{equation}

For another such hyperbolic metric space $(\Xcst, \dst)$ with base point $\ost$,
let $\qst$ denote the quasi-metric in $\partial \Xcst$.
In the product space,
for $\xib_i=(\xi_i, \eta_i) \in \(\Xc\cup\partial \Xc\) \times \(\Xcst\cup\partial \Xcst\)$ and $i=1, 2$,
let
\[
\qb(\xib_1, \xib_2):=\max\big\{q(\xi_1, \xi_2), \qst(\eta_1, \eta_2)\big\}.
\]
By \eqref{Eq:q} extended to $\Xc\cup \partial \Xc$ and $\Xcst \cup \partial \Xcst$, there exists a constant $C:=C_{\qb}>0$ such that
\[
\qb(\xb, \yb) \le C\max\{\qb(\xb, \zb), \qb(\zb, \yb)\} \quad \text{for all $\xb, \yb, \zb \in \(\Xc\cup\partial \Xc\)\times\(\Xcst \cup \partial \Xcst\)$}.
\]
Let $\Bb(\xib, r)$ denote the ball in $\(\Xc\cup\partial \Xc\) \times \(\Xcst\cup\partial \Xcst\)$ with respect to $\qb$.

\subsection{Hausdorff dimensions}

Let $(\Mcc, q)$ be a compact metrizable space $\Mcc$ with a quasi-metric $q$.
It is basically intended as $(\partial \Xc, q)$ or $(\partial \Xc\times \partial \Xcst, \qb)$.
For a set $E$ in $\Mcc$, 
let $\dim E$ denote the Hausdorff dimension of $E$ with respect to the quasi-metric $q$.
The definition is recalled briefly.
Let 
$|E|:=\sup\{q(\xi, \eta) \ : \ \xi, \eta \in E\}$.
For all $\alpha, \Delta \in \R_+$ with $\Delta>0$,
let
\[
\Hc^\alpha_\Delta(E):=\inf\Big\{\sum_{i=0}^\infty |E_i|^\alpha \ : \ E \subset \bigcup_{i=0}^\infty E_i \ \text{and} \ |E_i| \le \Delta\Big\}.
\]
The {\bf $\alpha$-dimensional Hausdorff measure} of a set $E$ is defined by
\[
\Hc^\alpha(E):=\lim_{\Delta \to 0}\Hc^\alpha_\Delta(E)=\sup_{\Delta>0}\Hc^\alpha_\Delta(E).
\]
Moreover the {\bf Hausdorff dimension} of a set $E$ is defined by
\[
\dim E:=\sup\Big\{\alpha \ge 0 \ : \ \Hc^\alpha(E)>0\Big\}=\inf\Big\{\alpha \ge 0 \ : \ \Hc^\alpha(E)=0\Big\}.
\]
Let $\nu$ be a Borel probability measure on $\partial \Xc$. 
The {\bf upper Hausdorff dimension} of $\nu$ is
\[
\overline \dim\, \nu:=\inf\Big\{\dim E \ : \ \text{$E$ is Borel and $\nu(\Mcc \setminus E)=0$} \Big\},
\]
and the {\bf lower Hausdorff dimension} of $\nu$ is
\[
\underline \dim\, \nu:=\inf\Big\{\dim E \ : \ \text{$E$ is Borel and $\nu(E)>0$}\Big\}.
\]
If the upper and lower Hausdorff dimensions of $\nu$ coincide,
then the common value is called the {\bf Hausdorff dimension} of $\nu$ and is denoted by $\dim\,\nu$.
The following is a fundamental lemma which relates pointwise behaviors of a measure to Hausdorff dimensions.
This is called the Billingsley lemma (in the case of Euclidean spaces).

\begin{lemma}[cf.\ Section 8.7 in \cite{Heinonen}]\label{Lem:Frostman}
For every Borel probability measure $\nu$ on $\Mcc$,
if for $\alpha_1, \alpha_2 \in \R_+$,
\[
\alpha_1 \le \liminf_{r \to 0}\frac{\log \nu(B(\xi, r))}{\log r} \le \alpha_2 \quad \text{for $\nu$-almost every $\xi \in \Mcc$},
\]
then $\alpha_1 \le \overline \dim\,\nu \le \alpha_2$.
\end{lemma}

It is deduced that
\[
\overline \dim\,\nu=\sup_{\text{$\nu$-a.e.\ $\xi$}}\liminf_{r \to 0}\frac{\log \nu(B(\xi, r))}{\log r}
\quad \text{and} \quad
\underline \dim\,\nu=\inf_{\text{$\nu$-a.e.\ $\xi$}}\liminf_{r \to 0}\frac{\log \nu(B(\xi, r))}{\log r}.
\]
In the above, $\sup_{\text{$\nu$-a.e.\ $\xi$}}$ and $\inf_{\text{$\nu$-a.e.\ $\xi$}}$ denote the essential supremum and the essential infimum relative to $\nu$ respectively.
A Borel probability measure $\nu$ on $\Mcc$ is {\bf exact dimensional} if the following limit
exists and is constant $\nu$-almost everywhere on $\Mcc$:
\[
\lim_{r\to 0}\frac{\log \nu(B(\xi, r))}{\log r}.
\]
In that case, the Hausdorff dimension of $\nu$ exists and equals the constant.

A metric space $(\Xc, d)$ is called {\bf bounded growth at some scale}
if there exist constants $r, R \in \R_+$ with $0<r<R$ and $N \in \Z_+$ such that every open ball of radius $R$ is covered by at most $N$ open balls of radius $r$.
The examples include Gromov hyperbolic Riemannian manifolds whose sectional curvature is uniformly bounded from below and from above and Cayley graphs of hyperbolic groups.

\begin{lemma}\label{Lem:BonkSchramm}
Let $(\Xc, d)$ and $(\Xcst, \dst)$ be hyperbolic metric spaces with bounded growth at some scale.
There exists a constant $L\ge 1$ such that the following holds
for every Borel probability measure $\nu$ on $\partial \Xc\times \partial \Xcst$ and for every Borel set $F$ in $\partial \Xc\times \partial \Xcst$ with $\nu(F)>0$.
For $\nu$-almost every $\xib \in F$, there exists a constant $r(\xib)>0$ such that for every $r \in (0, r(\xib))$,
\[
\nu(F \cap \Bb(\xib, L r)) \ge \frac{9}{10}\nu(\Bb(\xib, r)).
\]
\end{lemma}

\proof
The assumption on $(\Xc, d)$ implies that for every $\alpha \in (0, 1)$ there exists a bi-Lipschitz embedding $f$ from $(\partial \Xc, q^\alpha)$ to some finite dimensional standard Euclidean space $({\tt \Eo}, \|\cdot\|_{{\tt \Eo}})$ (cf.\ \cite[Theorem 9.2]{BonkSchramm} and \cite[2.6.\ Proposition]{Assouad}).
More precisely, there exists a constant $L^0 \ge 1$ such that for all $\xi, \eta \in \partial \Xc$,
\[
(1/L^0)q(\xi, \eta)^{\alpha} \le \|f(\xi)-f(\eta)\|_{{\tt E}_0} \le L^0 q(\xi, \eta)^{\alpha}.
\]
Similarly, for $(\Xcst, \dst)$ there exists a bi-Lipschitz embedding $\fst$ from $(\partial \Xcst, \qst^\alpha)$ into some Euclidean space $({\tt \Est}, \|\cdot\|_{{\tt \Est}})$ with a Lipschitz constant $\Lst \ge 1$.
Let 
\[
{\bm f}:\partial \Xc\times \partial \Xcst \to {\tt E}:={\tt \Eo} \times {\tt \Est}, \quad (\xi, \eta)\mapsto(f(\xi), \fst(\eta)).
\]
The map ${\bm f}$ is a homeomorphism onto its image since $\partial \Xc \times \partial \Xcst$ is compact.
The product space ${\tt E}$ is endowed with the maximum norm $\|\cdot\|_{{\tt E}}$ of the factors.
Let $B_{{\tt E}}(\vb, r)$ denote the ball in ${\tt E}$ with respect to the norm.
It holds that
\begin{equation}\label{Eq:L}
(1/L)\qb(\xib, \etab)^{\alpha} \le \|{\bm f}(\xib)-{\bm f}(\etab)\|_{{\tt E}} \le L \qb(\xib, \etab)^{\alpha} \quad \text{for $\xib, \etab \in \partial \Xc\times \partial \Xcst$},
\end{equation}
where $L:=\max\{L^0, \Lst\}$.
The pushforward ${\bm f}_\ast \nu$ satisfies that ${\bm f}_\ast \nu(B_{{\tt E}}({\bm f}(\xib), r))>0$ for all $r>0$ and for $\nu$-almost all $\xib \in F$. This follows since $\nu(F)>0$ and the intersection of $F$ and the support of $\nu$ has a positive $\nu$-measure, it holds that $\nu(\Bb(\xib, r))>0$ for all $r>0$ and for $\nu$-almost every $\xib \in F$.
The Lebesgue differentiation theorem on ${\bm f}_\ast \nu$ yields
\[
\lim_{r \to 0}\frac{{\bm f}_\ast \nu({\bm f}(F)\cap B_{{\tt E}}({\bm f}(\xib),r))}{{\bm f}_\ast \nu(B_{{\tt E}}({\bm f}(\xib),r))}=1 \quad \text{for $\nu$-almost every $\xib \in F$}.
\]
By \eqref{Eq:L}, it holds that
\[
\liminf_{r \to 0}\frac{\nu(F\cap B(\xib, (Lr)^{1/{\alpha}}))}{\nu(B(\xib, (r/L)^{1/{\alpha}}))}\ge 1 \quad \text{for $\nu$-almost every $\xib \in F$}.
\]
Hence for $\nu$-almost every $\xib \in F$ there exists some $r(\xib)>0$ such that
\[
\nu(F\cap B(\xib, L^{2/\alpha}r))\ge \frac{9}{10}\nu(B(\xib, r)) \quad \text{for all $r \in (0, r(\xib))$}.
\]
Shifting the constant $L^{2/\alpha}$ to $L$ deduces the claim.
\qed

\subsection{Random walks}\label{Sec:RW}

Let $\Gammab$ be a countable group.
Further let $\Omega:=\Gammab^{\Z_+}$ be the product space endowed with the $\sigma$-algebra $\Fc$ generated by cylinder sets.
For a probability measure $\pi$ on $\Gammab$,
let $\pi^{\Z_+}$ be the product measure on $\Gammab^{\Z_+}$.
Let us define the map $w:\Gammab^{\Z_+} \to \Omega$, $\{\xb_n\}_{n \in \Z_+}\mapsto \{\wb_n\}_{n \in \Z_+}$ where $\wb_0:=\id$ (the identity element) and 
\[
\wb_n:=\xb_1\cdots \xb_n \quad \text{for $n=1, 2, \dots$}.
\]
The pushforward of $\pi^{\Z_+}$ by the map $w$ is denoted by $\Pb$.
The probability space $(\Omega, \Fc, \Pb)$ is a standard probability space; this is the most basic space in the following discussion.
The maps $\Omega \to \Gammab$, $\{\wb_n\}_{n \in \Z_+} \mapsto \wb_n$ defines a Markov chain $\{\wb_n\}_{n \in \Z_+}$ called a $\pi$-{\bf random walk} starting from $\id$.

For a hyperbolic metric space $(\Xc, d)$, let $\Isom \Xc$ denote the isometry group.
A probability measure $\mu$ (with a countable support) on $\Isom \Xc$ is called {\bf non-elementary} if the group generated by the support of $\mu$ (as a group) contains a free group of rank $2$.

Let $(\Xc, d)$ and $(\Xcst, \dst)$ be hyperbolic metric spaces, and $\Gammab$ be a countable subgroup of $\Isom \Xc \times \Isom \Xcst$.
Further let $\pi$ be a probability measure on $\Gammab$ 
such that the pushforwards $\mu$ and $\must$ by the projections from
$\Isom \Xc \times \Isom \Xcst$ to $\Isom \Xc$ and to $\Isom \Xcst$ respectively are non-elementary.
In this setting, a $\pi$-random walk $\{\wb_n\}_{n \in \Z_+}$ starting from $\id$ yields by letting $\wb_n=(w_n, \wst_n)$,
a $\mu$-random walk $\{w_n\}_{n \in \Z_+}$ with $w_0=\id$ and a $\must$-random walk $\{\wst_n\}_{n \in \Z_+}$ with $\wst_0=\id$.
For fixed base points $o \in \Xc$ and $\ost \in \Xcst$,
let
\[
\zb_n:=(z_n, \zst_n), \quad \text{where $z_n:=w_n\cdot o$ and $\zst_n:=\wst_n \cdot \ost$ for $\wb_n=(w_n, \wst_n)$}.
\]
The assumption that $\mu$ and $\must$ are non-elementary implies that $\Pb$-almost surely there exist $z_\infty \in \partial \Xc$ and $\zst_\infty \in \partial \Xcst$ such that
$z_n \to z_\infty$ in $\Xc\cup \partial \Xc$ and $\zst_n \to \zst_\infty$ in $\Xcst \cup \partial \Xcst$
as $n \to \infty$ respectively.
Let $\nu_\pi$ be the distribution of $(z_\infty, \zst_\infty)$ on $\partial \Xc \times \partial \Xcst$, and $\nu_\mu$ and $\nu_{\must}$ be the distributions of $z_\infty$ and of $\zst_\infty$ respectively.
The probability measure $\nu_\pi$ is called the {\bf harmonic measure} for the $\pi$-random walk.
Similarly, $\nu_\mu$ and $\nu_{\must}$ are called the harmonic measures for the $\mu$-random walk and for the $\must$-random walk respectively.
Let us denote the measurable map by
\[
\bnd=(\rmbnd, \rmbndst):\Omega \to \partial \Xc\times \partial \Xcst, \quad \wb\mapsto \zb_\infty:=(z_\infty, \zst_\infty).
\]
The harmonic measures $\nu_\pi$, $\nu_\mu$ and $\nu_{\must}$ are obtained as the pushforwards of $\Pb$ by $\bnd$, by $\rmbnd$ and by $\rmbndst$ respectively. 

The group $\Isom \Xc \times \Isom \Xcst$ acts on $\partial \Xc \times \partial \Xcst$ coordinatewise
and $\Gammab$ acts on $\partial \Xc \times \partial \Xcst$.
It satisfies that $\nu_\pi$ is {\bf $\pi$-stationary}, i.e.,
\[
\nu_\pi=\sum_{\xb \in \Gammab}\pi(\xb)\xb\nu_\pi \quad \text{where $\xb\nu_\pi=\nu_\pi \circ \xb^{-1}$}.
\]
Similarly, $\Gammab$ acts on $\partial \Xc$ and on $\partial \Xcst$ through the projection from $\Isom \Xc \times \Isom \Xcst$ to each one of factors,
and thus $\nu_\mu$ and $\nu_{\must}$ are also $\pi$-stationary, i.e.,
\begin{equation}\label{Eq:pi-stationary_factors}
\nu_\mu=\sum_{\xb \in \Gammab}\pi(\xb)\xb\nu_\mu \quad \text{and} \quad 
\nu_{\must}=\sum_{\xb \in \Gammab}\pi(\xb)\xb\nu_{\must}.
\end{equation}
Since $\mu$ and $\must$ are marginals of $\pi$, these further lead to
\[
\nu_\mu=\sum_{x \in \supp \mu}\mu(x)x\nu_\mu \quad \text{and} \quad \nu_{\must}=\sum_{\xst \in \supp \must}\must(\xst)\xst\nu_{\must}.
\]

Let us define the metric on $\Xc\times\Xcst$ by 
\[
\db(\zb_1, \zb_2):=\max\big\{d(z_1, z_2), \dst(\zst_1, \zst_2)\big\}
\quad \text{for $\zb_i=(z_i, \zst_i) \in \Xc\times \Xcst$ and $i=1, 2$}.
\]
A probability measure $\pi$ on $\Gammab < \Isom \Xc \times \Isom \Xcst$ has a {\bf finite first moment} if
\[
\sum_{\xb \in \Gammab}\db(\ob, \xb\cdot\ob)\pi(\xb)<\infty, \quad \text{where $\ob=(o, \ost)$}.
\]
This condition is independent of the choice of the point $\ob=(o, \ost)$.
Let us assume that $\pi$ has a finite first moment.
The Kingman subadditive ergodic theorem implies that 
the following limits exist and are constant $\Pb$-almost everywhere:
\[
l(\Xc, \mu):=\lim_{n \to \infty}\frac{1}{n}d(o, z_n) \quad \text{and} \quad l(\Xcst, \must):=\lim_{n \to \infty}\frac{1}{n}\dst(\ost, \zst_n).
\]
The limits $l(\Xc, \mu)$ and $l(\Xcst, \must)$ are called the {\bf drifts} of $\{z_n\}_{n \in \Z_+}$ and $\{\zst_n\}_{n \in \Z_+}$ respectively.
In the case when $\mu$ and $\must$ are non-elementary,
then $l(\Xc, \mu)>0$ and $l(\Xcst, \must)>0$ (\cite[Theorem 7.3]{Kaimanovich-hyperbolic}, and see also \cite[Theorem 1.1]{Gouezel_Exp} for a more recent account).

\subsection{Conditional processes and their entropies}

If $\pi$ has a finite first moment,
then the entropy $H(\pi):=-\sum_{\xb \in \Gammab}\pi(\xb)\log \pi(\xb)$ is finite \cite[Section VII, B]{Derriennic}.
The Shannon theorem for random walks says that 
for such $\pi$, the following limit exists and is constant $\Pb$-almost everywhere:
\begin{equation}\label{Eq:entropy}
h(\pi):=\lim_{n \to \infty}-\frac{1}{n}\log \pi_n(\wb_n).
\end{equation}
See \cite[Theorem 2.1]{KaimanovichVershik} and \cite[Section IV]{Derriennic80}.
The limit $h(\pi)$ is called the {\bf asymptotic entropy} for $\pi$-random walk.
Let $h(\mu)$ and $h(\must)$ be the asymptotic entropies for $\mu$-random walk and $\must$-random walk respectively; they exist and are defined in the same way.
We will also use a conditional version of the notion.
First we introduce a conditional process and then define the conditional entropy.

Recall that $(\Omega, \Fc, \Pb)$ is a standard probability space.
Let $\sigma(\bnd)$ be the $\sigma$-algebra generated by the measurable map $\bnd: \Omega \to \partial \Xc \times \partial \Xcst$.
Disintegrating the measure $\Pb$ with respect to $\sigma(\bnd)$ yields the system of conditional probability measures $\{\Pb^{\bnd(\wb)}\}_{\wb \in \Omega}$.
More precisely,
for every $A \in \Fc$, 
the map $\Omega \to \R$, $\wb \mapsto \Pb^{\bnd(\wb)}(A)$ is $\sigma(\bnd)$-measurable, and
\[
\Pb=\int_\Omega \Pb^{\bnd(\wb)}\,d\Pb(\wb).
\]
Noting that $\nu_\pi=\bnd_\ast \Pb$,
let us write $\{\Pb^{\xi, \eta}\}_{(\xi, \eta) \in \partial \Xc \times \partial \Xcst}$ and
\[
\Pb=\int_{\partial \Xc \times \partial \Xcst}\Pb^{\xi, \eta}\,d\nu_\pi(\xi, \eta).
\]
Similarly, disintegrating $\Pb$ with respect to the $\sigma$-algebra $\sigma(\rmbndst)$ generated by $\rmbndst$ yields the system of conditional probability measures $\{\Pb^\eta\}_{\eta \in \partial \Xcst}$ satisfying that
\begin{equation}\label{Eq:disintegration_eta}
\Pb=\int_{\partial \Xcst}\Pb^\eta\,d\nu_{\must}(\eta).
\end{equation}
Let us disintegrate the harmonic measure $\nu_\pi$.
In the present setting, $\partial \Xc\times \partial \Xcst$ is a compact metrizable space and thus $\partial \Xc \times \partial \Xcst$ endowed with the Borel $\sigma$-algebra is a standard Borel space.
The probability measure $\nu_\pi$ is disintegrated with respect to the $\sigma$-algebra generated by the projection $\partial \Xc \times \partial \Xcst \to \partial \Xcst$.
This yields the system of conditional probability measures $\{\nu_\pi^\eta\}_{\eta \in \partial \Xcst}$ such that
\[
\nu_\pi=\int_{\partial \Xcst}\nu_\pi^\eta\,d\nu_{\must}(\eta).
\]
Moreover, it satisfies that 
$\nu_\pi^\eta(\partial \Xc\times \partial \Xcst)=\nu_\pi^\eta(\partial \Xc\times \{\eta\})=1$
for $\nu_{\must}$-almost every $\eta \in \partial \Xcst$.
These disintegrations lead to by the Fubini theorem,
\begin{equation}\label{Eq:P_xi_eta}
\Pb^\eta=\int_{\partial \Xc \times \{\eta\}}\Pb^{\xi, \eta}\,d\nu_\pi^\eta(\xi) \quad \text{for $\nu_{\must}$-almost every $\eta \in \partial \Xcst$}.
\end{equation}

For $\nu_{\must}$-almost every $\eta \in \partial \Xcst$,
the conditional probability measure $\Pb^\eta$ coincides with the distribution of a conditional process on $\Gammab$.
This is a Markov chain whose transition probability is defined 
for $\nu_{\must}$-almost every $\eta \in \partial \Xcst$,
\[
p^\eta(\xb, \yb):=\pi(\xb^{-1}\yb)\frac{d\yb\nu_{\must}}{d\xb\nu_{\must}}(\eta) \quad \text{if $\pi(\xb^{-1}\yb)>0$, and $0$ if otherwise, for $\xb, \yb \in \Gammab$}.
\]
Note that $\yb\nu_{\must}$ is absolutely continuous with respect to $\xb\nu_{\must}$ and $d\yb\nu_{\must}/d\xb\nu_{\must}$ is well-defined $\nu_{\must}$-almost everywhere if $\pi(\xb^{-1}\yb)>0$
by \eqref{Eq:pi-stationary_factors}.
Moreover since 
\[
\frac{d\yb\nu_{\must}}{d\xb\nu_{\must}}(\eta)=\frac{d\xb^{-1}\yb\nu_{\must}}{d\nu_{\must}}(\xb^{-1}\eta) \quad \text{for $\nu_{\must}$-almost every $\eta \in \partial \Xcst$},
\]
the above $p^\eta(\xb, \yb)$ indeed defines a transition probability by the $\pi$-stationarity of $\nu_\must$.
Let
$\pi^\eta_n(\xb):=\Pb^\eta(\wb_n=\xb)$
for $\xb \in \Gammab$ and $n \in \Z_+$.
It holds that for $\nu_{\must}$-almost every $\eta \in \partial \Xcst$,
for every cylinder set $[\wb_0, \dots, \wb_n]$ in $(\Omega, \Fc, \Pb)$, 
\begin{equation}\label{Eq:Doob}
\Pb^\eta([\wb_0, \dots, \wb_n])=\Pb([\wb_0, \dots, \wb_n])\frac{d\wb_n \nu_\must}{d\nu_\must}(\eta).
\end{equation}
In particular, for $\nu_\must$-almost every $\eta \in \partial \Xcst$,
\[
\pi^\eta_n(\xb)=\pi_n(\xb)\frac{d\xb\nu_{\must}}{d\nu_{\must}}(\eta) \quad \text{for $\xb \in \Gammab$ and $n \in \Z_+$}.
\]
For more details, see \cite[Sections 3 and 4]{Kaimanovich-hyperbolic}.
There it is shown (in a more general setting) that
the Shannon theorem holds for the conditional process.
Namely,
for $\nu_{\must}$-almost every $\eta \in \partial \Xcst$,
the following limit exists and is constant $\Pb^\eta$-almost everywhere:
\begin{equation}\label{Eq:Shannon_condition}
h(\Pb^\eta):=\lim_{n \to \infty}-\frac{1}{n}\log \pi^\eta_n(\wb_n).
\end{equation}
Furthermore, the limit is obtained as
\[
h(\Pb^\eta)=h(\pi)-\sum_{\xb \in \Gammab}\pi(\xb)\int_{\partial \Xcst}\log \frac{d\xb\nu_{\must}}{d\nu_{\must}}(\eta)\,d\xb\nu_{\must}(\eta).
\]
See \cite[Theorem 4.5]{Kaimanovich-hyperbolic}.
The {\bf differential entropy} for the pair $(\partial \Xcst, \must)$ is defined by
\[
h(\partial \Xcst, \must):=\sum_{x \in \supp \must}\must(x)\int_{\partial \Xcst}\log \frac{dx\nu_{\must}}{d\nu_{\must}}(\eta)\,dx\nu_{\must}(\eta).
\]
Since $\must$ is a marginal of $\pi$,
it holds that
\[
h(\Pb^\eta)=h(\pi)-h(\partial \Xcst, \must) \quad \text{for $\nu_{\must}$-almost every $\eta \in \partial \Xcst$}.
\]
Let us mention that the differential entropy arises in the theory of Poisson boundary in the following way:
It has been proven that 
$h(\partial \Xcst, \must) \le h(\must)$
and the equality holds if and only if $(\partial \Xcst, \nu_{\must})$ is a Poisson boundary for the $\must$-random walk \cite[Theorem 4.6]{Kaimanovich-hyperbolic}.
In the present setting, since $(\partial \Xcst, \nu_{\must})$ is $\pi$-stationary \eqref{Eq:pi-stationary_factors}, it holds that
$h(\partial \Xcst, \must)=h(\partial \Xcst, \pi)$, and that $h(\partial \Xcst, \pi)\le h(\pi)$, 
where the equality holds if and only if $(\partial \Xcst, \nu_\must)$ is a Poisson boundary for $(\Gammab, \pi)$.

\section{Exact dimension of conditional measures}\label{Sec:dimension}

For a proper $C$-roughly geodesic hyperbolic metric space $(\Xc, d)$ for some $C \in \R_+$ with a fixed base point $o \in \Xc$ and for a $\mu$-random walk $\{w_n\}_{n \in \Z_+}$,
the following {\bf ray approximation} holds for $z_n=w_n\cdot o$:
If $\mu$ is non-elementary and has a finite first moment,
then $\Pb$-almost surely there exists a $C$-rough geodesic ray $\gamma_{z_\infty}$ such that for $l:=l(\Xc, \mu)$ of $\{z_n\}_{n \in \Z_+}$,
\begin{equation}\label{Eq:ray}
d(z_n, \gamma_{z_\infty}(ln))=o(n).
\end{equation}
See \cite[Theorem 7.3]{Kaimanovich-hyperbolic}.
In fact, such an assignment $\xi\mapsto \gamma_\xi$ from $\partial \Xc$ to the space of $C$-rough geodesic rays from $o$ in $(\Xc, d)$ equipped with the topology of convergence on compact sets is chosen to be Borel measurable by the Borel selection theorem (cf.\ \cite[Section 3.2]{Tdim}).
In the same way, 
there is a Borel measurable map $\eta \mapsto \gamma_\eta$ from $\partial \Xcst$ to the space of $C$-rough geodesic rays from $\ost$ in $(\Xcst, \dst)$.
For the drift $\lst:=l(\Xcst, \must)$ of $\{\zst_n\}_{n \in \Z_+}$,
it holds that $\Pb$-almost surely,
\begin{equation}\label{Eq:ray_star}
\dst(\zst_n, \gamma_{\zst_\infty}(\lst n))=o(n).
\end{equation}

In the following subsections, for brevity, let
\[
\wbar h:=h(\pi)-h(\partial \Xcst, \must), \quad l:=l(\Xc, \mu) \quad  \text{and} \quad \lst:=l(\Xcst, \must).
\]

\subsection{Upper bounds on dimensions of conditional measures}

\begin{lemma}\label{Lem:dim_upper}
Let $\Gammab$ be a countable subgroup in $\Isom \Xc \times \Isom \Xcst$, and $\pi$ be a probability measure on $\Gammab$ with finite first moment and non-elementary marginals $\mu$ and $\must$. 
It holds that for $\nu_{\must}$-almost every $\eta \in \partial \Xcst$ and for $\nu_\pi^\eta$-almost every $\xib \in \partial \Xc\times \partial \Xcst$,
\[
\limsup_{r \to 0}\frac{\log \nu_\pi^\eta(\Bb(\xib, r))}{\log r} \le \frac{h(\pi)-h(\partial \Xcst, \must)}{l(\Xc, \mu)}.
\]
\end{lemma}

\proof
Recall that for $\nu_{\must}$-almost every $\eta \in \partial \Xcst$,
the distribution $\Pb^\eta$ is obtained by a Markov chain whose law at time $n$ is $\pi^\eta_n$ for $n \in \Z_+$.
For every $\e>0$ and every interval $I$ in $\Z_+$,
let
\[
A_{\e, I}:=\bigcap_{n \in I}\Big\{\wb \in \Omega \ : \ (z_n|z_{n+1})_o\ge (l-\e)n, \ \pi^{\rmbndst(\wb)}_n(\wb_n) \ge \exp(-n(\wbar h+\e))\Big\}.
\]
Note that $(z_n|z_{n+1})_o/n \to l$ as $n \to \infty$ almost surely in $\Pb$ since 
$\mu$ has a finite first moment and the $\mu$-random walk has the drift $l>0$.
By disintegration \eqref{Eq:disintegration_eta}, this together with \eqref{Eq:Shannon_condition} implies that
for every $\e\in (0, l)$ and for $\nu_{\must}$-almost every $\eta \in \partial \Xcst$,
\[
\Pb^\eta\Big(\bigcup_{N \in \Z_+}A_{\e, [N, \infty)}\Big)=1.
\]
Hence for every $\e>0$ and
for $\nu_{\must}$-almost every $\eta \in \partial \Xcst$,
there exists an $N_{\e, \eta} \in \Z_+$ such that
\[
\Pb^\eta\big(A_{\e, [N_{\e, \eta}, \infty)}\big) \ge 1-\e.
\]
Let $N:=N_{\e, \eta}$ and $A:=A_{\e, N_{\e, \eta}}$.
Further for all $n>N$,
let
\[
A_{[N, n)}:=A_{\e, [N, n)} \quad \text{and} \quad A_{[n, \infty)}:=A_{\e, [n, \infty)}.
\]

Note that $A=A_{[N, n)}\cap A_{[n, \infty)}$.
For $n \in \Z_+$ and $\wb \in \Omega$, let
\[
C_n(\wb):=\big\{\wb' \in \Omega \ : \ \wb_n'=\wb_n\big\}.
\]
This is the event where the position of the chain is $\wb_n$ at time $n$.
Since the conditional process is a Markov chain,
for $\nu_{\must}$-almost every $\eta \in \partial \Xcst$,
for all $\wb \in \Omega$ and all $n>N$,
\begin{equation}\label{Eq:Markov}
\Pb^\eta(A \mid C_n(\wb))=\Pb^\eta(A_{[N, n)}\mid C_n(\wb))\cdot \Pb^\eta(A_{[n, \infty)}\mid C_n(\wb)).
\end{equation}
Furthermore for $\nu_{\must}$-almost every $\eta \in \partial \Xcst$,
\begin{equation}\label{Eq:ANn}
\Pb^\eta(A_{[N, n)}\mid C_n(\wb))=\Pb^\eta(A_{[N, n)}\mid \sigma(\wb_n, \wb_{n+1},\dots)) \quad \text{almost everywhere in $\Pb^\eta$},
\end{equation}
where $\sigma(\wb_n, \wb_{n+1}, \dots)$ is the $\sigma$-algebra generated by $\wb_n, \wb_{n+1}, \dots$.
Similarly, one has
\begin{equation}\label{Eq:Aninfty}
\Pb^\eta(A_{[n, \infty)}\mid C_n(\wb))=\Pb^\eta(A_{[n, \infty)}\mid \sigma(\wb_0, \dots, \wb_n)) \quad \text{almost everywhere in $\Pb^\eta$},
\end{equation}
where $\sigma(\wb_0, \dots, \wb_n)$ is the $\sigma$-algebra generated by $\wb_0, \dots, \wb_n$.
Let us denote the tail $\sigma$-algebra by
\[
\Tc:=\bigcap_{n \in \Z_+}\sigma(\wb_n, \wb_{n+1}, \dots).
\]
Note that $\Pb^\eta$-almost everywhere on $A=A_{[N, n)}\cap A_{[n, \infty)}$,
\[
\Pb^\eta(A_{[N, n)}\mid \sigma(\wb_n, \wb_{n+1}, \dots))=\Pb^\eta(A\mid \sigma(\wb_n, \wb_{n+1}, \dots)).
\]
By \eqref{Eq:ANn}, the L\'evy downward theorem applied to the right hand side above shows that
for $\nu_{\must}$-almost every $\eta \in \partial \Xcst$,
\begin{equation}\label{Eq:ANn1}
\lim_{n \to \infty}\Pb^\eta(A_{[N, n)}\mid C_n(\wb))=\Pb^\eta(A\mid \Tc) \quad \text{almost everywhere in $\Pb^\eta$ on $A$}.
\end{equation}
Analogously, note that $\Pb^\eta$-almost everywhere on $A=A_{[N, n)}\cap A_{[n, \infty)}$,
\[
\Pb^\eta(A_{[n, \infty)}\mid \sigma(\wb_0, \dots, \wb_n))=\Pb^\eta(A\mid \sigma(\wb_0, \dots, \wb_n)).
\]
By \eqref{Eq:Aninfty},
the L\'evy upward theorem shows that
for $\nu_{\must}$-almost every $\eta \in \partial \Xcst$,
\begin{equation}\label{Eq:Aninfty1}
\lim_{n \to \infty}\Pb^\eta(A_{[n, \infty)}\mid C_n(\wb))=\1_A \quad \text{almost everywhere in $\Pb^\eta$ on $A$}.
\end{equation}
For $\nu_{\must}$-almost every $\eta \in \partial \Xcst$,
it holds that 
\begin{equation}\label{Eq:positive}
\Pb^\eta(A\mid \Tc)>0 \quad \text{almost everywhere in $\Pb^\eta$ on $A$}. 
\end{equation}
Indeed, if we define $\Nc:=\{\wb \in \Omega \ : \ \Pb^\eta(A\mid \Tc)=0\}$,
then $\Nc$ is $\Tc$-measurable and $\Pb^\eta(A\cap \Nc \mid \Tc)=0$, implying that $\Pb^\eta(A\cap \Nc)=0$.
Thus we have \eqref{Eq:positive}.
Therefore by \eqref{Eq:Markov}, \eqref{Eq:ANn1}, \eqref{Eq:Aninfty1} and \eqref{Eq:positive},
for $\nu_{\must}$-almost every $\eta \in \partial \Xcst$, 
\begin{equation}\label{Eq:main_limit}
\lim_{n \to \infty}\Pb^\eta(A\mid C_n(\wb))=\Pb^\eta(A\mid \Tc)>0 \quad \text{almost everywhere in $\Pb^\eta$ on $A$}.
\end{equation}

For every $\wb \in A$, it holds that
\[
q(z_n, z_\infty)=\exp(-(z_n|z_\infty)_o )\le Ce^{-(l-\e)n} \quad \text{for all $n>N$},
\]
where $C$ is independent of $\wb$.
This follows since $q$ with a power is bi-Lipschitz to a genuine metric (cf.\ Section \ref{Sec:hyperbolic}).
Further for $n>N$ and for $\Pb^\eta$-almost every $\wb' \in A\cap C_n(\wb)$,
\[
\wb_n'=\wb_n \quad \text{and} \quad q(z_n',z_\infty')\le Ce^{-(l-\e)n},
\]
where $z'_n:=w_n'\cdot o$ for $\wb_n=(w_n', {\wst_n}^{'})$ and $z_n' \to z_\infty'$ as $n \to \infty$ in $\Xc \cup\partial \Xc$.
Since $\wb_n'=\wb_n$ and thus $z_n'=z_n$, 
by \eqref{Eq:delta} extended on $\Xc\cup \partial \Xc$,
\[
q(z_\infty, z_\infty')\le C e^{\delta-(l-\e)n}, \quad \text{i.e.,} \quad z_\infty' \in B(z_\infty, Ce^{\delta-(l-\e)n}) \quad \text{for all $n>N$}.
\]
Hence for $\nu_{\must}$-almost every $\eta \in \partial \Xcst$,
\[
\Pb^\eta(A\cap C_n(\wb)) \le \nu_\pi^\eta(B(z_\infty, Ce^{\delta-(l-\e)n})\times \partial \Xcst) \quad \text{for all $n>N$}.
\]
The right hand side coincides with $\nu_\pi^\eta(\Bb(\zb_\infty, Ce^{\delta-(l-\e)n}))$ where $\zb_\infty=(z_\infty, \zst_\infty)$
since $\nu_\pi^\eta$ is supported in $\partial \Xc \times \{\eta\}$ and $\eta=\zst_\infty$ for $\nu_{\must}$-almost every $\eta \in \partial \Xcst$.
For $\Pb^\eta$-almost every $\wb \in A$,
it holds that for all $n>N$,
\[
\Pb^\eta(A\cap C_n(\wb))=\Pb^\eta(A\mid C_n(\wb))\cdot \Pb^\eta(C_n(\wb)) \ge \Pb^\eta(A\mid C_n(\wb))\cdot e^{-n(\wbar h+\e)}.
\]
This shows that for $\nu_{\must}$-almost every $\eta \in \partial \Xcst$, for $\Pb^\eta$-almost every $\wb \in A$ and for all $n>N$,
\[
\nu_\pi^\eta(\Bb(\zb_\infty, Ce^{\delta-(l-\e)n})) \ge \Pb^\eta(A\mid C_n(\wb))\cdot e^{-n(\wbar h+\e)}.
\]
By \eqref{Eq:main_limit}, for $\nu_{\must}$-almost every $\eta \in \partial \Xcst$ and for $\Pb^\eta$-almost every $\wb \in A$,
\begin{equation}\label{Eq:limsup_A}
\limsup_{r \to 0}\frac{\log \nu_\pi^\eta(\Bb(\zb_\infty, r))}{\log r} \le \frac{\wbar h+\e}{l-\e}.
\end{equation}
This follows first for the sequence $r_n\to 0$ as $n \to \infty$ where $r_n:=Ce^{\delta-(l-\e)n}$ and then for $r>0$ and $r\to 0$ by noting that $r_{n+1}=e^{-(l-\e)}r_n$.
Recall that $\Pb^\eta(A) \ge 1-\e$ for the event $A=A_{\e, [N_{\e, \eta}, \infty)}$ for $\nu_{\must}$-almost every $\eta \in \partial \Xcst$.
For an arbitrary decreasing sequence $\e_n \to 0$ as $n \to \infty$,
one has
$\Pb^\eta(\bigcap_{m \in \Z_+}\bigcup_{n \ge m}A_{\e_n, [N_{\e_n, \eta}, \infty)})=1$
for $\nu_{\must}$-almost every $\eta \in \partial \Xcst$.
Therefore by \eqref{Eq:limsup_A} for $\nu_{\must}$-almost every $\eta \in \partial \Xcst$, 
\[
\limsup_{r\to 0}\frac{\log \nu_\pi^\eta(\Bb(\zb_\infty, r))}{\log r}\le \frac{\wbar h}{l} \quad \text{almost everywhere in $\Pb^\eta$}.
\]
Noting that for $\nu_{\must}$-almost every $\eta \in \partial \Xcst$, the distribution of $\zb_\infty$ is $\nu_\pi^\eta$,
we obtain for $\nu_{\must}$-almost every $\eta \in \partial \Xcst$,
\[
\limsup_{r \to 0}\frac{\log \nu_\pi^\eta(\Bb(\xib, r))}{\log r} \le \frac{\wbar h}{l} \quad \text{for $\nu_\pi^\eta$-almost every $\xib \in \partial \Xc\times \partial \Xcst$}.
\]
This concludes the claim.
\qed

We use the following version of Lemma \ref{Lem:dim_upper} in Section \ref{Sec:product}.

\begin{lemma}\label{Lem:upper_F}
In the same setting as in Lemma \ref{Lem:dim_upper}, 
if for $\nu_{\mu^\star}$-almost every $\eta \in \partial \Xcst$
there exists a Borel set $F_\eta$ in $\partial \Xc\times \partial \Xcst$ such that $\nu_\pi^\eta(F_\eta)>0$,
then for $\nu_{\mu^\star}$-almost every $\eta \in \partial \Xcst$ and for $\nu_\pi^\eta$-almost every $\xib \in F_\eta$,
\[
\limsup_{r \to 0}\frac{\log \nu_\pi^\eta(\Bb(\xib, r)\cap F_\eta)}{\log r} \le \frac{h(\pi)-h(\partial \Xcst, \mu^\star)}{l(\Xc, \mu)}.
\]
\end{lemma}

\proof
This follows from Lemmas \ref{Lem:BonkSchramm} and \ref{Lem:dim_upper}. 
\qed

\subsection{Lower bounds on dimensions of conditional measures}

\begin{lemma}\label{Lem:lower_strong}
Let $\Gammab$ be a countable subgroup in $\Isom \Xc \times \Isom \Xcst$ with finite exponential growth relative to $(\Xc \times \Xcst, \db)$, and $\pi$ be a probability measure on $\Gammab$ with finite first moment and non-elementary marginals $\mu$ and $\must$.
For every $\e>0$, there exist 
\begin{enumerate}
\item an $N \in \Z_+$, 
\item a Borel set $D$ in $\partial \Xcst$ with $\nu_\must(D) \ge 1-\e$,
and 
\item a Borel set $F$ in $\partial \Xc \times \partial \Xcst$ with $\nu_\pi^\eta(F) \ge 1-\e$ for $\nu_\must$-almost every $\eta \in D$ and $\nu_\pi(F) \ge 1-\e$,
\end{enumerate}
such that the following holds:
For $\nu_\must$-almost every $\eta \in D$, 
for all $\xi \in \partial \Xc$ and all $n \ge N$,
\[
\nu_\pi^\eta\(B(\xi, e^{-l n})\times \partial \Xcst \cap F\) \le C_\e e^{-n(\wbar h-\e)},
\]
where $C_\e$ is a constant depending only on $\e$.
\end{lemma}

\proof
For every $\e>0$ and every $N \in \Z_+$,
let
\[
A_{\e, N}:=\bigcap_{n \ge N}\Big\{\wb \in \Omega \ : \ \db(\zb_n, (\gamma_{z_\infty}(l n), \gamma_{\zst_\infty}(\lst n))) \le \e n, \ \pi^{\rmbndst(\wb)}_n(\wb_n)\le \exp(-n(\wbar h-\e))\Big\}.
\]
The disintegration formula \eqref{Eq:disintegration_eta} implies that
for all event $A \subset \Omega$ if $\Pb(A)=1$, then $\Pb^\eta(A)=1$ for $\nu_{\must}$-almost every $\eta \in \partial \Xcst$.
This together with \eqref{Eq:ray} and \eqref{Eq:ray_star},
and further \eqref{Eq:Shannon_condition}
imply that for every $\e>0$,
\[
\Pb^\eta\Big(\bigcup_{N \in \Z_+}A_{\e, N}\Big)=1 \quad \text{for $\nu_{\must}$-almost every $\eta \in \partial \Xcst$}.
\]
For $N \in \Z_+$, let
\[
D_{\e, N}:=\Big\{\eta \in \partial \Xcst \ : \ \Pb^\eta(A_{\e, N}) \ge 1-\e\Big\}.
\]
Since $A_{\e, N}$ is increasing and $\Pb^\eta(A_{\e, N}) \to 1$ monotonically as $N \to \infty$ for $\nu_\must$-almost every $\eta \in \partial \Xcst$,
there exists an $N_\e \in \Z_+$ such that
\begin{equation}\label{Eq:Lem:lower_strong_D}
\nu_\must\(D_{\e, N_\e}\) \ge 1-\e.
\end{equation}
Let $N:=N_\e$, $D:=D_{\e, N_\e}$ and $A:=A_{\e, N_\e}$.
It holds that
\begin{equation}\label{Eq:P^eta(A)}
\Pb^\eta(A) \ge 1-\e \quad \text{for $\nu_\must$-almost every $\eta \in D$}.
\end{equation}
Let us define
\[
F_\e:=\Big\{(\xi, \eta) \in \partial \Xc \times \partial \Xcst \ : \ \Pb^{\xi, \eta}(A) \ge \e \Big\},
\]
and $F:=F_\e$.
We claim that 
\begin{equation}\label{Eq:nu^eta}
\nu_\pi^\eta(F)\ge 1-2 \e \quad \text{for $\nu_{\must}$-almost every $\eta \in D$}.
\end{equation}
Indeed, by \eqref{Eq:P^eta(A)} and by \eqref{Eq:P_xi_eta},
for $\nu_{\must}$-almost every $\eta \in D$,
\[
1-\e \le \Pb^\eta(A)=\int_{\partial \Xc\times \{\eta\}}\Pb^{\xi, \eta}(A)\,d\nu_\pi^\eta(\xi)=\int_F \Pb^{\xi, \eta}(A)\,d\nu_\pi^\eta(\xi)+\int_{\partial \Xc\times \partial \Xcst\setminus F}\Pb^{\xi, \eta}(A)\,d\nu_\pi^\eta(\xi).
\]
This implies that
\[
1-\e \le \nu_\pi^\eta(F)+\e \cdot \nu_\pi^\eta(\partial \Xc\times \partial \Xcst \setminus F)\le \nu_\pi^\eta(F)+\e,
\]
showing \eqref{Eq:nu^eta}.

Furthermore it holds that
\begin{equation}\label{Eq:Lem:lower_strong_F}
\nu_\pi(F) \ge 1-3\e.
\end{equation}
This follows since $\nu_\pi^\eta(F) \ge (1-2\e)\1_D$ by \eqref{Eq:nu^eta},
integration with respect to $\nu_\must$ yields
\[
\nu_\pi(F)=\int_{\partial \Xcst}\nu_\pi^\eta(F)\,d\nu_\must(\eta) \ge (1-2\e)\nu_\must(D) \ge (1-2\e)(1-\e) \ge 1-3\e,
\]
where the second inequality uses \eqref{Eq:Lem:lower_strong_D}.
It holds that for $\nu_{\must}$-almost every $\eta \in \partial \Xcst$,
for every $\xi \in \partial \Xc$ and $R>0$, 
\begin{align}\label{Eq:main_A}
&\nu_\pi^\eta(\Oc(\gamma_\xi(l n), R)\times \partial \Xcst\cap F)
=\Pb^\eta(\{\zb_\infty \in \Oc(\gamma_\xi(ln), R)\times \partial \Xcst \cap F\})\nonumber\\
&\qquad \qquad=\Pb^\eta(\{\zb_\infty \in \Oc(\gamma_\xi(l n), R)\times \partial \Xcst\cap F\}\cap A)\nonumber\\
&\qquad \qquad \qquad \qquad \qquad+\Pb^\eta(\{\zb_\infty \in \Oc(\gamma_\xi(l n), R)\times \partial \Xcst\cap F\}\cap A^{\sf c}).
\end{align}
In the above, $A^{\sf c}$ denote the complement event of $A$, and $\zb_\infty=(z_\infty, \zst_\infty)$.
First let us bound the first term in \eqref{Eq:main_A}.
On the event $A$, 
if $z_\infty \in \Oc(\gamma_\xi(l n), R)$ and $n \ge N$,
then
\[
d(z_n, \gamma_\xi(l n))\le d(z_n, \gamma_{z_\infty}(l n))+d(\gamma_{z_\infty}(l n), \gamma_\xi(l n)) \le \e n +C_R,
\]
where $C_R:=2R+2C$ since $\gamma_\xi$ is a $C$-rough geodesic ray.
Moreover, since $\nu_{\must}$-almost every $\eta \in \partial \Xcst$ it holds that $\zst_\infty=\eta$ almost surely in $\Pb^\eta$,
it holds that $\Pb^\eta$-almost everywhere on $A$ for all $n \ge N$,
\[
\dst(\zst_n, \gamma_\eta(\lst n)) \le \e n.
\]
Hence $\Pb^\eta$-almost everywhere on $A$ for all $n \ge N$,
\[
\db(\zb_n, \gamma_{\xi, \eta, n}) \le \e n + C_R \quad \text{where $\gamma_{\xi, \eta, n}:=(\gamma_\xi(l n), \gamma_\eta(\lst n))$}.
\]
This shows that letting $\Bb(\xb, r)$ denote the ball in $(\Xc \times \Xcst, \db)$,
we have for all $n \ge N$,
\begin{align*}
&\Pb^\eta(\{\zb_\infty \in \Oc(\gamma_\xi(l n), R)\times \partial \Xcst \cap F\}\cap A)\\
&\qquad \qquad \qquad\le \Pb^\eta(\{\zb_n \in \Bb(\gamma_{\xi, \eta, n}, \e n +C_R)\}\cap\{\pi^\eta_n(\wb_n)\le \exp(-n(\wbar h-\e))\}).
\end{align*}
The right hand side is at most 
$\sum \pi^\eta_n(\xb)$ where the summation runs over all $\xb \in \Gammab$ such that $\xb \cdot \ob \in \Bb(\gamma_{\xi, \eta, n}, \e n +C_R)$ and $\pi^\eta_n(\xb) \le \exp(-n(\wbar h-\e))$.
This is at most
\[
\#\big\{\xb \in \Gammab \ : \ \xb \cdot \ob \in \Bb(\gamma_{\xi, \eta, n}, \e n +C_R)\big\}\cdot e^{-n(\wbar h-\e)}
\le c e^{c(\e n+C_R)}\cdot e^{-n(\wbar h-\e)},
\]
for a constant $c>0$ since $\Gammab$ has a finite exponential growth relative to $(\Xc \times \Xcst, \db)$.
Thus for $\nu_{\must}$-almost every $\eta \in \partial \Xcst$ and for every $\xi \in \partial \Xc$, for all $n \ge N$,
\begin{equation}\label{Eq:main_A1}
\Pb^\eta(\{\zb_\infty \in \Oc(\gamma_\xi(l n), R)\times \partial \Xcst \cap F\}\cap A) \le c e^{c(\e n +C_R)}\cdot e^{-n(\wbar h-\e)}.
\end{equation}

Next let us bound the second term in \eqref{Eq:main_A}.
By \eqref{Eq:P_xi_eta}, 
it holds that for $\nu_\must$-almost every $\eta \in D$,
\begin{align}\label{Eq:main_A2}
&\Pb^\eta(\{\zb_\infty \in \Oc(\gamma_\xi(l n), R)\times \partial \Xcst \cap F\}\cap A^{\sf c})\nonumber\\
&=\int_{\Oc(\gamma_\xi(l n), R)\times \partial \Xcst \cap F}\Pb^{\zeta, \eta}(A^{\sf c})\,d\nu_\pi^\eta(\zeta)
\le (1-\e)\cdot \nu_\pi^\eta(\Oc(\gamma_\xi(l n), R)\times \partial \Xcst \cap F).
\end{align}
In the above, the inequality holds since $\Pb^{\zeta, \eta}(A^{\sf c}) \le 1-\e$ for $\nu_\pi^\eta$-almost every $(\zeta, \eta) \in F$ for $\nu_\must$-almost every $\eta \in D$ by the definition of $F$.

Finally, combining \eqref{Eq:main_A}, \eqref{Eq:main_A1} and \eqref{Eq:main_A2} yields
for $\nu_\must$-almost every $\eta \in D$, for every $\xi \in \partial \Xc$ and
for all $n \ge N$,
\[
\nu_\pi^\eta(\Oc(\gamma_\xi(l n), R)\times \partial \Xcst \cap F)
\le c e^{c(\e n+C_R)}\cdot e^{-n(\wbar h-\e)}+(1-\e)\cdot \nu_\pi^\eta(\Oc(\gamma_\xi(l n), R)\times \partial \Xcst\cap F).
\]
Therefore for $\nu_{\must}$-almost every $\eta \in D$ and for every $\xi \in \partial \Xc$, for all $n \ge N$,
\begin{equation}\label{Eq:main_A3}
\e \cdot \nu_\pi^\eta\(\Oc(\gamma_\xi(l n), R)\times \partial \Xcst\cap F\)
\le c e^{c(\e n+C_R)}\cdot e^{-n(\wbar h-\e)}.
\end{equation}
Note that $B(\xi, C^{-1}e^{- l n+R})\cap \partial \Xc\subset \Oc(\gamma_\xi(l n), R)$ by \eqref{Eq:shadows_balls},
where we choose a large enough constant $R$ so that $e^R/C\ge 1$ and $C$ depends only on the hyperbolicity constant.
In \eqref{Eq:main_A3}, the constant $c$ depends only on $\Gammab$ and $(\Xc\times \Xcst, \db)$.
For every $\e>0$, we argue with $\e'= \e/(3+c)$.
By \eqref{Eq:Lem:lower_strong_D}, \eqref{Eq:nu^eta} and \eqref{Eq:Lem:lower_strong_F}, 
we obtain 
\[
\nu_\must(D) \ge 1-\e'\ge 1-\e, \quad \nu_\pi^\eta(F)\ge 1-2\e' \ge 1-\e \quad \text{for $\nu_\must$-almost every $\eta \in D$},
\]
and $\nu_\pi(F) \ge 1-3\e'\ge 1-\e$.
Further by \eqref{Eq:main_A3},
for $\nu_\must$-almost every $\eta \in D$, for every $\xi \in \partial \Xc$ and for all $n \ge N$,
\[
\nu_\pi^\eta\(B(\xi, e^{-ln})\cap F\)
\le (c/\e') e^{cC_R}\cdot e^{-n(\wbar h-\e)}.
\]
Defining the constant $C_\e:=(c/\e') e^{cC_R}$ yields the claim.
\qed

\begin{lemma}\label{Lem:dim_lower}
In the same setting as in Lemma \ref{Lem:lower_strong},
it holds that for $\nu_{\must}$-almost every $\eta \in \partial \Xcst$ and for $\nu_\pi^\eta$-almost every $\xib \in \partial \Xc\times \partial \Xcst$,
\[
\liminf_{r \to 0}\frac{\log \nu_\pi^\eta(\Bb(\xib, r))}{\log r} \ge \frac{h(\pi)-h(\partial \Xcst, \must)}{l(\Xc, \mu)}.
\]
\end{lemma}
\proof
By Lemma \ref{Lem:lower_strong},
for every $\e>0$ there exist a Borel set $D$ in $\partial \Xcst$ with $\nu_\must(D) \ge 1-\e$ and a Borel set $F$ in $\partial \Xc\times \partial \Xcst$ with $\nu_\pi(F) \ge 1-\e$ such that the following holds:
For $\nu_{\must}$-almost every $\eta \in D$ and for every $\xi\in \partial \Xc$,
\[
\liminf_{n \to \infty}\frac{\log \nu_\pi^\eta(\Bb((\xi, \eta), e^{-ln})\cap F)}{-ln}\ge \frac{\wbar h-\e}{l}.
\]
In fact, in the above the sequence $r_n:=e^{-ln}$ for $n \in \Z_+$ is replaced by positive reals $r$ tending to $0$ since $r_{n+1}=e^{-l}r_n$ for all $n \in \Z_+$.
Applying Lemma \ref{Lem:BonkSchramm} to the measures $\nu_\pi^\eta$ and $F$ implies the following:
There exists a constant $L\ge 1$ such that for $\nu_{\must}$-almost every $\eta \in D$, for $\nu_\pi^\eta$-almost every $(\xi, \eta)\in F$ and for a constant $r(\xi, \eta)>0$,
\[
\nu_\pi^\eta(\Bb((\xi, \eta), L r)\cap F) \ge \frac{9}{10}\nu_\pi^\eta(\Bb((\xi, \eta), r)) \quad \text{for all $r \in (0, r(\xi, \eta))$}.
\]
Hence for $\nu_{\must}$-almost every $\eta \in D$ and for $\nu_\pi^\eta$-almost every $\xib=(\xi, \eta)\in F$,
\[
\liminf_{r\to 0}\frac{\log \nu_\pi^\eta(\Bb(\xib, r))}{\log r}\ge \frac{\wbar h-\e}{l}.
\]
Since for every $\e>0$ there exists such an $F$ denoted by $F_\e$ with $\nu_\pi^\eta(F_\e) \ge 1-\e$,
for $\nu_{\must}$-almost every $\eta \in D$,
it holds that $\nu_\pi^\eta(\bigcap_{m\in \Z_+}\bigcup_{n \ge m}F_{\e_n})=1$ for an arbitrary decreasing sequence $\e_n \to 0$ as $n \to \infty$.
Therefore it follows that for $\nu_{\must}$-almost every $\eta \in D$ and for $\nu_\pi^\eta$-almost every $\xib \in \partial \Xc\times \partial \Xcst$,
\[
\liminf_{r \to 0}\frac{\log \nu_\pi^\eta(\Bb(\xib, r))}{\log r}\ge \frac{\wbar h}{l}.
\]
For every $\e>0$ there exists such a $D$ with $\nu_\must(D) \ge 1-\e$,
the above holds for $\nu_\must$-almost every $\eta \in \partial \Xcst$.
This concludes the claim.
\qed

\subsection{Proofs of Theorems \ref{Thm:exact_intro} and \ref{Thm:exact}}

\proof[Proof of Theorem \ref{Thm:exact}]
Since $\mu$ and $\must$ are non-elementary and have finite first moments,
the drift $l(\Xc, \mu)$ is finite and positive, further the asymptotic entropy $h(\pi)$ and the differential entropy $h(\partial \Xcst, \must)$ are finite.
The first claim follows from Lemmas \ref{Lem:dim_upper} and \ref{Lem:dim_lower}.
The second claim follows from Lemma \ref{Lem:Frostman}.
\qed

\proof[Proof of Theorem \ref{Thm:exact_intro}]
Since a probability measure $\must$ on $\Gammast$ is non-elementary with finite first moment, it holds that $h(\must)=h(\partial \Gammast, \must)$ since $(\partial \Gammast, \nu_\must)$ is a Poisson boundary for $(\Gammast, \must)$ \cite[Theorem 7.4]{Kaimanovich-hyperbolic}.
The claim follows from Theorem \ref{Thm:exact} by
applying $\Gamma$ and $\Gammast$ endowed with left invariant hyperbolic metrics quasi-isometric to word metrics respectively to $(\Xc, d)$ and $(\Xcst, \dst)$.
\qed

\section{Exact dimension of harmonic measures in products spaces}\label{Sec:product}

As in Section \ref{Sec:dimension}, for brevity, let
\[
l:=l(\Xc, \mu), \quad \lst:=l(\Xcst, \must), \quad \hst:=h(\partial \Xcst, \must) \quad \text{and} \quad \wbar h:=h(\pi)-\hst.
\]

\subsection{Upper bounds on dimensions of harmonic measures in product spaces}

The proof of the following proposition is inspired by \cite[Section 8]{LedrappierLessaExact}.

\begin{proposition}\label{Prop:upper_product}
Let $\Gammab$ be a countable subgroup in $\Isom \Xc\times \Isom \Xcst$, and $\pi$ be a probability measure on $\Gammab$ with finite first moment and non-elementary marginals $\mu$ and $\must$.
If $l(\Xc, \mu)\ge l(\Xcst, \must)$,
then
it holds that
for $\nu_\pi$-almost every $\xib \in \partial \Xc \times \partial \Xcst$,
\[
\limsup_{r \to 0}\frac{\log \nu_\pi\(\Bb(\xib, r)\)}{\log r} \le \frac{h(\pi)-h(\partial \Xcst, \must)}{l(\Xc, \mu)}+\frac{h(\partial \Xcst, \must)}{l(\Xcst, \must)}.
\]
\end{proposition}

We assume that $l\ge \lst$: if otherwise we argue after exchanging the notations $l$ and $\lst$.
Fix an arbitrary $\e\in(0, \lst)$.
Let 
\[
r_n:=e^{-(\lst-\e)n} \quad \text{for $n \in \Z_+$}.
\]
Let us define (recalling that $\zb_t=\wb_t\cdot \ob$)
\[
A_{\e, n}:=\bigcap_{t\ge n}\Big\{\wb \in \Omega \ : \ \text{$\zb_\infty=(z_\infty, \zst_\infty)$ exists and $\qb(\zb_t, \zb_\infty) \le r_n$}\Big\}.
\]
\begin{lemma}\label{Lem:A}
The events $A_{\e, n}$ are increasing in $n \in \Z_+$, and
it holds that 
\begin{equation*}\label{Eq:Lem:A}
\Pb\Big(\bigcup_{n \in \Z_+}A_{\e, n}\Big)=1.
\end{equation*}
\end{lemma}
\proof
By definition $A_{\e, n}$ are increasing in $n \in \Z_+$.
For $\Pb$-almost every $\wb \in \Omega$, for all large enough $t \in \Z_+$ (recalling that $z_t=w_t\cdot o$ and $\zst_t=\wst_t\cdot \ost$),
\[
(z_t|z_\infty)_o\ge (l-\e) t \quad \text{and} \quad (\zst_t|\zst_\infty)_{\ost}\ge (\lst -\e)t.
\]
In which case, $\max \{e^{-(z_t|z_\infty)_o}, e^{-(\zst_t|\zst_\infty)_{\ost}}\} \le e^{-(\lst-\e)t}$ since $l \ge \lst$, showing the claim.
\qed

\medskip

For each $\eta \in \partial \Xcst$
and $n \in \Z_+$,
let
\[
E_{\e, n}(\eta):=\Big\{\wb \in \Omega \ : \ \Pb^\eta([\wb_0, \dots, \wb_n]\cap A_{\e, n})\le e^{n(\hst+\e)}\Pb([\wb_0, \dots, \wb_n]\cap A_{\e, n})\Big\},
\]
and $E_{\e, [n, \infty)}(\eta):=\bigcap_{t \ge n}E_{\e, t}(\eta)$.

\begin{lemma}\label{Lem:Eneta}
For each $\eta \in \partial \Xcst$, the events $E_{\e, [n, \infty)}(\eta)$ are increasing in $n \in \Z_+$, and for $\nu_\must$-almost every $\eta \in \partial \Xcst$,
\begin{equation*}\label{Eq:Eneta}
\Pb^\eta\Big(\bigcup_{n \in \Z_+}E_{\e, [n, \infty)}(\eta)\Big)=1.
\end{equation*}
\end{lemma}

\proof
In the following, let $A_n:=A_{\e, n}$ for $n \in \Z_+$.
By \eqref{Eq:Doob},
for $\nu_{\must}$-almost every $\eta \in \partial \Xcst$, for every cylinder set $[\wb_0, \dots, \wb_n]$ in $(\Omega, \Fc, \Pb)$,
\[
\Pb^\eta([\wb_0, \dots, \wb_n])=\Pb([\wb_0, \dots, \wb_n])\frac{d\wb_n \nu_{\must}}{d\nu_{\must}}(\eta).
\]
Note that for $\Pb$-almost every $\wb \in \Omega$, for all $n \in \Z_+$, 
\[
\Pb([\wb_0, \dots, \wb_n])>0 \quad \text{and} \quad \Pb^{\rmbndst(\wb)}([\wb_0, \dots, \wb_n])>0.
\]
The Birkhoff ergodic theorem implies that $\Pb$-almost every $\wb \in \Omega$,
\[
\lim_{n \to \infty}\frac{1}{n}\log \frac{d\wb_n\nu_{\must}}{d\nu_{\must}}(\rmbndst(\wb)) = \hst.
\]
See \cite[the proof of Theorem 4.5]{Kaimanovich-hyperbolic}.
This implies that by disintegration of $\Pb$ into $\Pb^\eta$ for $\eta \in \partial \Xcst$, for $\nu_\must$-almost every $\eta \in \partial \Xcst$, for $\Pb^\eta$-almost every $\wb \in \Omega$,
\[
\lim_{n \to \infty}\frac{1}{n}\log \frac{d\wb_n\nu_\must}{d\nu_\must}(\eta)=\hst.
\]
Therefore for $\nu_\must$-almost every $\eta \in \partial \Xcst$, for $\Pb^\eta$-almost every $\wb \in \Omega$,
\begin{equation}\label{Eq:Lem:Eneta0}
\lim_{n \to \infty}\frac{1}{n}\log \frac{\Pb^{\eta}([\wb_0, \dots, \wb_n])}{\Pb([\wb_0, \dots, \wb_n])}=\hst.
\end{equation}

Fix an arbitrary $N \in \Z_+$. 
For all $n\in [N, \infty)\cap\Z_+$, for all cylinder set $[\wb_0, \dots, \wb_n]$ in $(\Omega, \Fc, \Pb)$
of positive $\Pb$-measure,
it holds that
\[
\frac{\Pb([\wb_0, \dots, \wb_n]\cap A_{N})}{\Pb([\wb_0, \dots, \wb_n])} \le \frac{\Pb([\wb_0, \dots, \wb_n]\cap A_{n})}{\Pb([\wb_0, \dots, \wb_n])} \le 1.
\]
The left most side equals $\Pb(A_{N}\mid \sigma(\wb_0, \dots, \wb_n))$ almost everywhere in $\Pb$.
The martingale convergence theorem yields
\[
\lim_{n \to \infty}\Pb(A_{N}\mid \sigma(\wb_0, \dots, \wb_n))=\1_{A_{N}} \quad \text{for $\Pb$-almost every $\wb \in \Omega$}.
\]
Since $N$ is arbitrary and $\Pb(\bigcup_{N \in \Z_+}A_{N})=1$,
it holds that
\begin{equation*}
\lim_{n \to \infty}\frac{\Pb([\wb_0, \dots, \wb_n]\cap A_{n})}{\Pb([\wb_0, \dots, \wb_n])}=\1 \quad \text{for $\Pb$-almost every $\wb \in \Omega$}.
\end{equation*}
By disintegration of $\Pb$ into $\Pb^\eta$ for $\eta \in \partial \Xcst$,
for $\nu_\must$-almost every $\eta \in \partial \Xcst$,
\begin{equation}\label{Eq:Lem:EnetaA}
\lim_{n \to \infty}\frac{\Pb([\wb_0, \dots, \wb_n]\cap A_{n})}{\Pb([\wb_0, \dots, \wb_n])}=\1 \quad \text{for $\Pb^\eta$-almost every $\wb \in \Omega$}.
\end{equation}

Applying the same discussion to $\Pb^\eta$ as for $\Pb$, we obtain for $\nu_\must$-almost every $\eta \in \partial \Xcst$,
\begin{equation}\label{Eq:Lem:EnetaAeta}
\lim_{n \to \infty}\frac{\Pb^\eta([\wb_0, \dots, \wb_n]\cap A_{n})}{\Pb^\eta([\wb_0, \dots, \wb_n])}=\1 \quad \text{for $\Pb^\eta$-almost every $\wb \in \Omega$}.
\end{equation}
Combining \eqref{Eq:Lem:Eneta0}, \eqref{Eq:Lem:EnetaA} and \eqref{Eq:Lem:EnetaAeta} yields
for $\nu_\must$-almost every $\eta \in \partial \Xcst$,
\begin{equation}\label{Eq:Lem:Eneta}
\lim_{n \to \infty}\frac{1}{n}\log \frac{\Pb^\eta([\wb_0, \dots, \wb_n]\cap A_n)}{\Pb([\wb_0, \dots, \wb_n]\cap A_n)}=\hst \quad \text{for $\Pb^\eta$-almost every $\wb \in \Omega$}.
\end{equation}
By definition for each $\eta \in \partial \Xcst$ the events $E_{\e, [n, \infty)}(\eta)$ are increasing in $n \in \Z_+$ respectively, and by \eqref{Eq:Lem:Eneta} for $\nu_\must$-almost every $\eta \in \partial \Xcst$, 
\[
\Pb^\eta\Big(\bigcup_{n \in \Z_+}E_{\e, [n, \infty)}(\eta)\Big)=1,
\]
as claimed.
\qed

\proof[Proof of Proposition \ref{Prop:upper_product}]
Fix an arbitrary $\e \in (0, \lst)$ and recall that $r_n:=e^{-(\lst-\e)n}$ for $n \in \Z_+$.
Let $A_n$, $E_n(\eta)$ and $E_{[n, \infty)}(\eta)$ denote $A_{\e, n}$, $E_{\e, n}(\eta)$ and $E_{\e, [n, \infty)}(\eta)$ respectively for brevity.
Lemma \ref{Lem:A} implies that by disintegration of $\Pb$ into $\Pb^\eta$ for $\eta \in \partial \Xcst$,
for $\nu_\must$-almost every $\eta \in \partial \Xcst$,
\[
\Pb^\eta\Big(\bigcup_{n \in \Z_+}A_{n}\Big)=1.
\]
By this together with Lemma \ref{Lem:Eneta},
for $\nu_\must$-almost every $\eta \in \partial \Xcst$,
there exists an $N_{\e, \eta} \in \Z_+$ such that
\begin{equation}\label{Eq:AE}
\Pb^\eta\(E_{[N_{\e, \eta}, \infty)}(\eta)\cap A_{N_{\e, \eta}}\) \ge 1-\e.
\end{equation}
Let $N:=N_{\e, \eta}$, $E_{[N, \infty)}(\eta):=E_{[N_{\e, \eta}, \infty)}(\eta)$ and $A_N:=A_{N_{\e, \eta}}$.
Further let
\[
F_\eta:=\Big\{\xib\in \partial \Xc\times \partial \Xcst \ : \ \Pb^{\xib}\(E_{[N, \infty)}(\eta)\cap A_N\) \ge \e \Big\}.
\]
Note that $F_\eta$ is a Borel measurable set in $\partial \Xc\times \partial \Xcst$ since for each $B \in \Fc$ in $(\Omega, \Fc, \Pb)$, the map $\xib \mapsto \Pb^{\xib}(B)$ is Borel measurable.
By \eqref{Eq:AE} and by disintegration of $\Pb^\eta$ into $\Pb^{\xib}$ for $\xib \in \partial \Xc\times \partial \Xcst$,
for $\nu_\must$-almost every $\eta \in \partial \Xcst$,
\[
1-\e \le \Pb^\eta\(E_{[N, \infty)}(\eta)\cap A_N\)=\int_{\partial \Xc\times \partial \Xcst}\Pb^{\xib}\(E_{[N, \infty)}(\eta)\cap A_N\)\,d\nu_\pi^\eta(\xib)
\le \nu_\pi^\eta(F_\eta)+\e \nu_\pi^\eta(F_\eta^{\sf c}).
\]
Therefore for $\nu_\must$-almost every $\eta \in \partial \Xcst$,
\begin{equation}\label{Eq:Feta}
\nu_\pi^\eta(F_\eta)\ge 1-2\e.
\end{equation}
Furthermore, for $\nu_\must$-almost every $\eta \in \partial \Xcst$ and for every $\xib_0 \in \partial \Xc\times \partial \Xcst$, 
\begin{align*}
\nu_\pi^\eta\(\Bb(\xib_0, r_n)\cap F_\eta\) 
&= \Pb^\eta\(\{\zb_\infty \in \Bb(\xib_0, r_n)\cap F_\eta\}\cap E_{[N, \infty)}(\eta)\cap A_N\)\\
&\qquad +\Pb^\eta\(\{\zb_\infty \in \Bb(\xib_0, r_n)\cap F_\eta\}\cap (E_{[N, \infty)}(\eta)\cap A_N)^{\sf c}\).
\end{align*}
By definition of $F_\eta$, it holds that
\begin{align*}
&\Pb^\eta\(\{\zb_\infty \in \Bb(\xib_0, r_n)\cap F_\eta\}\cap (E_{[N, \infty)}(\eta)\cap A_N)^{\sf c}\)\\
&=
\int_{\Bb(\xib_0, r_n)\cap F_\eta}\Pb^{\xib}\((E_{[N, \infty)}(\eta)\cap A_N)^{\sf c}\)\,d\nu_\pi^\eta(\xib)
\le (1-\e)\nu_\pi^\eta\(\Bb(\xib_0, r_n)\cap F_\eta\).
\end{align*}
In summary, for all $n \ge N$,
\begin{align*}
\e \nu_\pi^\eta\(\Bb(\xib_0, r_n)\cap F_\eta\) 
\le \Pb^\eta\(\{\zb_\infty \in \Bb(\xib_0, r_n)\cap F_\eta\}\cap E_{[N, \infty)}(\eta)\cap A_N\).
\end{align*}
By Lemmas \ref{Lem:A} and \ref{Lem:Eneta}, the events $A_n$ and $E_{[n, \infty)}(\eta)$ are increasing in $n \in \Z_+$, and $E_{[n, \infty)}(\eta) \subset E_n(\eta)$ by the definition.
Thus, for $\nu_\must$-almost every $\eta \in \partial \Xcst$, for every $\xib_0 \in \partial \Xc\times \partial \Xcst$, and for all $n \ge N=N_{\e, \eta}$,
\begin{equation}\label{Eq:upper_AE0}
\e \nu_\pi^\eta\(\Bb(\xib_0, r_n)\cap F_\eta\) 
\le \Pb^\eta\(\{\zb_\infty \in \Bb(\xib_0, r_n)\}\cap E_{n}(\eta)\cap A_n\).
\end{equation}
By definition of $A_n$, if $A_n$ holds, then $\zb_n \in \Bb(\zb_\infty, r_n)$, whence for $C:=C_{\qb}>0$,
\begin{equation}\label{Eq:upper_AE1}
\Pb^\eta\(\{\zb_\infty \in \Bb(\xib_0, r_n)\}\cap E_{n}(\eta)\cap A_n\)
\le \Pb^\eta\(\{\zb_n \in \Bb(\xib_0, Cr_n)\}\cap E_{n}(\eta)\cap A_n\).
\end{equation}
For every $\xib_0 \in \partial \Xc\times \partial \Xcst$ and for every $n \in \Z_+$,
let
\[
B_n(\xib_0):=\big\{\zb_n \in \Bb(\xib_0, Cr_n)\big\}.
\]
For each such fixed $\xib_0$ and $n$, the event $B_n(\xib_0)$ is $\sigma(\wb_0, \dots, \wb_n)$-measurable.
Further each fixed $\eta$ and $n$, the event $E_n(\eta)$ is $\sigma(\wb_0, \dots, \wb_n)$-measurable.
Hence for each fixed $\xib_0$, $\eta$ and $n$, the event $B_n(\xib_0)\cap E_n(\eta)$ is $\sigma(\wb_0, \dots, \wb_n)$-measurable and is obtained as a (countable) sum of cylinder sets $[\wb_0, \dots, \wb_n]$ in $(\Omega, \Fc, \Pb)$.
Decomposing the event into a sum of cylinder sets yields
\begin{align*}
\Pb^\eta\(B_n(\xib_0)\cap E_n(\eta)\cap A_n\)
=\sum_{[\wb_0, \dots, \wb_n] \subset B_n(\xib_0)\cap E_n(\eta)}\Pb^\eta\([\wb_0, \dots, \wb_n]\cap A_n\).
\end{align*}
In the right hand side,
each summand is at most $e^{n(\hst+\e)}\Pb\([\wb_0, \dots, \wb_n]\cap A_n\)$
since $[\wb_0, \dots, \wb_n] \subset E_n(\eta)$.
Furthermore since $\Pb\([\wb_0, \dots, \wb_n]\cap A_n\)$ over those cylinder sets add up to $\Pb(B_n(\xib_0)\cap E_n(\eta)\cap A_n)$,
it holds that
\begin{equation}\label{Eq:upper_AE2}
\Pb^\eta\(B_n(\xib_0)\cap E_n(\eta)\cap A_n\)\le e^{n(\hst+\e)}\Pb\(B_n(\xib_0)\cap E_n(\eta)\cap A_n\).
\end{equation}
By the definitions of $B_n(\xib_0)$ and $A_n$, 
if $B_n(\xib_0)\cap A_n$ holds, then $\zb_\infty \in \Bb(\xib_0, C^2r_n)$.
This in particular implies that
\begin{equation}\label{Eq:upper_AE3}
\Pb\(B_n(\xib_0)\cap E_n(\eta)\cap A_n\)\le \Pb\(\{\zb_\infty \in \Bb(\xib_0, C^2 r_n)\}\)=\nu_\pi\(\Bb(\xib_0, C^2 r_n)\)
\end{equation}

Combining \eqref{Eq:upper_AE0}, \eqref{Eq:upper_AE1}, \eqref{Eq:upper_AE2} and \eqref{Eq:upper_AE3} implies that 
for $\nu_\must$-almost every $\eta \in \partial \Xcst$, for every $\xib_0 \in \partial \Xc\times \partial \Xcst$, and for all $n \ge N=N_{\e, \eta}$,
\begin{equation}\label{Eq:upper0}
\e \nu_\pi^\eta\(\Bb(\xib_0, r_n)\cap F_\eta\) \le e^{n(\hst+\e)}\nu_\pi\(\Bb(\xib_0, C^2 r_n)\).
\end{equation}

Recall that $\wbar h=h(\pi)-\hst$.
By Lemma \ref{Lem:upper_F},
for $\nu_\must$-almost every $\eta \in \partial \Xcst$,
for $\nu_\pi^\eta$-almost every $\xib_0 \in F_\eta$, there exists an $N_{\eta, \xib_0} \in \Z_+$ such that
for all $n \ge N_{\eta, \xib_0}$,
\[
\(\frac{\wbar h}{l}+\e\)\log r_n \le \log \nu_\pi^\eta\(\Bb(\xib_0, r_n)\cap F_\eta\).
\]

This together with \eqref{Eq:upper0} shows that 
for $\nu_\must$-almost every $\eta \in \partial \Xcst$,
for $\nu_\pi^\eta$-almost every $\xib_0 \in F_\eta$
(recalling that $r_n=e^{-(\lst-\e)n}$),
\[
\limsup_{n \to \infty}\frac{\log \nu_\pi\(\Bb(\xib_0, C^2 r_n)\)}{\log r_n}
\le \frac{\wbar h}{l}+\e+\frac{\hst+\e}{\lst-\e}.
\]
Recall that $\nu_\pi^\eta(F_\eta)\ge 1-2\e$ by \eqref{Eq:Feta} and that this holds for arbitrary $\e\in (0, \lst)$.
Therefore after replacing the sequence $r_n$ by reals $r$ tending to $0$,
for $\nu_\must$-almost every $\eta \in \partial \Xcst$,
for $\nu_\pi^\eta$-almost every $\xib \in \partial \Xc\times \partial \Xcst$,
\[
\limsup_{r \to 0}\frac{\log \nu_\pi\(\Bb\(\xib, r\)\)}{\log r}\le \frac{\wbar h}{l}+\frac{\hst}{\lst}.
\]
The disintegration of $\nu_\pi$ into $\nu_\pi^\eta$ shows that the above holds for $\nu_\pi$-almost every $\xib \in \partial \Xc\times \partial \Xcst$, concluding the claim.
\qed

\subsection{Lower bounds on dimensions of harmonic measures in product spaces}

\begin{proposition}\label{Prop:lower_product}
Let $\Gammab=\Gamma \times \Gammast$ where $\Gamma$ and $\Gammast$ are countable subgroups in $\Isom \Xc$ and in $\Isom \Xcst$ with finite exponential growth relative to $(\Xc, d)$ and to $(\Xcst, \dst)$ respectively.
For every probability measure $\pi$ on $\Gammab$ with finite first moment and non-elementary marginals $\mu$ and $\must$, 
it holds that for $\nu_\pi$-almost every $\xib \in \partial \Xc \times \Xcst$,
\[
\liminf_{r \to 0}\frac{\log \nu_\pi\(\Bb(\xib, r)\)}{\log r} \ge \frac{h(\pi)-h(\must)}{l(\Xc, \mu)}+\frac{h(\must)}{l(\Xcst, \must)}.
\]
\end{proposition}

\proof
Note that $\Gammab$ has a finite exponential growth relative to $(\Xc \times \Xcst, \db)$ by the assumption.
By Lemma \ref{Lem:lower_strong},
for every $\e>0$ there exist an $N \in \Z_+$, a Borel set $D$ in $\partial \Xcst$ with $\nu_\must(D) \ge 1-\e$ and a Borel set $F$ in $\partial \Xc \times \partial \Xcst$ with $\nu_\pi(F) \ge 1-\e$ as stated.
If we define $\hat F:=(\partial \Xc \times D) \cap F$,
then
\begin{equation}\label{Eq:hatF}
\nu_\pi(\hat F) \ge 1-2\e.
\end{equation}
This follows since $\nu_\pi(\partial \Xc \times D)=\nu_\must(D) \ge 1-\e$ and $\nu_\pi(F) \ge 1-\e$.

Let $r_n:=e^{-l n}$ for $n \in \Z_+$.
For every $\xib=(\xi, \xist) \in \partial \Xc \times \partial \Xcst$,
by disintegration of $\nu_\pi$ into $\nu_\pi^\eta$ for $\eta \in \partial \Xcst$,
\[
\nu_\pi\big(B(\xi, r_n)\times B(\xist, r_n)\cap \hat F\big)
=\int_{B(\xist, r_n)\cap D}\nu_\pi^\eta\(B(\xi, r_n) \times \partial \Xcst \cap F\)\,d\nu_\must(\eta).
\]
By Lemma \ref{Lem:lower_strong}, 
for $\nu_\must$-almost every $\eta \in D$, for every $\xi \in \partial \Xc$ and for all $n \ge N$,
\[
\nu_\pi^\eta\(B(\xi, r_n)\times \partial \Xcst \cap F\) \le C_\e e^{-n(\wbar h-\e)}.
\]
Therefore
for every $\xib=(\xi, \xist) \in \partial \Xc\times \partial \Xcst$ and for all $n \ge N$,
\begin{equation}\label{Eq:lower0}
\nu_\pi\big(B(\xi, r_n)\times B(\xist, r_n)\cap \hat F\big) \le C_\e e^{-n(\wbar h-\e)} \nu_\must\(B(\xist, r_n)\cap D\).
\end{equation}
By the dimension formula for $(\Gammast, \must)$ in \cite[Theorem 1.2]{Tdim},
for $\nu_\must$-almost every $\xist \in \partial \Xcst$,
\[
\lim_{n \to \infty}\frac{\log \nu_\pi\(B(\xist, r_n)\)}{\log r_n} = \frac{h(\must)}{\lst}.
\]
(The proof presented there is for geodesic spaces, but it is adapted to roughly geodesic spaces $\Xc$.)
Since $\nu_\pi\(B(\xist, r_n)\cap D\) \le \nu_\pi\(B(\xist, r_n)\)$ and $\log r_n<0$,
for $\nu_\must$-almost every $\xist \in \partial \Xcst$,
\[
\liminf_{n \to \infty}\frac{\log \nu_\pi\(B(\xist, r_n)\cap D\)}{\log r_n} \ge \frac{h(\must)}{\lst}.
\]
This together with \eqref{Eq:lower0} implies that
since $\nu_\pi(\partial \Xc \times D)=\nu_\must(D)$ and $\Bb(\xib, r_n)=B(\xi, r_n)\times B(\xist, r_n)$,
for $\nu_\pi$-almost every $\xib=(\xi, \xist) \in \partial \Xc\times D$ (recalling that $r_n=e^{-ln}$),
\begin{equation}\label{Eq:lower1}
\liminf_{n \to \infty}\frac{\log \nu_\pi\big(\Bb(\xib, r_n)\cap \hat F\big)}{\log r_n}\ge \frac{\wbar h-\e}{l}+\frac{h(\must)}{\lst}.
\end{equation}
By Lemma \ref{Lem:BonkSchramm}, for $\nu_\pi$-almost every $\xib \in \hat F$ (where $\hat F \subset \partial \Xc \times D$),
\begin{equation}\label{Eq:lower2}
\liminf_{n \to \infty}\frac{\log \nu_\pi\big(\Bb(\xib, r_n)\big)}{\log r_n}\ge \frac{\wbar h-\e}{l}+\frac{h(\must)}{\lst}.
\end{equation}
By \eqref{Eq:hatF}, one has $\nu_\pi(\hat F) \ge 1-2\e$, and for every $\e>0$ there exists such an $\hat F$ in $\partial \Xc \times \partial \Xcst$.
Thus after replacing the sequence $r_n$ for $n \in \Z_+$ by positive reals $r$ tending to $0$, we obtain
for $\nu_\pi$-almost every $\xib \in \partial \Xc \times \partial \Xcst$,
\[
\liminf_{r \to 0}\frac{\log \nu_\pi\big(\Bb(\xib, r)\big)}{\log r}\ge \frac{\wbar h}{l}+\frac{h(\must)}{\lst}.
\]
This concludes the claim.
\qed

\subsection{Exact dimension and the proof of Theorem \ref{Thm:product_intro}}

\begin{theorem}\label{Thm:product}
Let $(\Xc, d)$ and $(\Xcst, \dst)$ be roughly geodesic hyperbolic metric spaces with bounded growth at some scale.
Let $\Gammab=\Gamma \times \Gammast$ where $\Gamma$ and $\Gammast$ are countable subgroups in $\Isom \Xc$ and in $\Isom \Xcst$ with finite exponential growth relative to $(\Xc, d)$ and to $(\Xcst, \dst)$ respectively.
For every probability measure $\pi$ on $\Gammab$ with finite first moment and non-elementary marginals $\mu$ and $\must$, 
the harmonic measure $\nu_\pi$ on $\partial \Xc \times \partial \Xcst$ is exact dimensional.
Moreover, if $l(\Xc,\mu)\ge l(\Xcst, \must)$,
then it holds that
\[
\dim \nu_\pi= \frac{h(\pi)-h(\must)}{l(\Xc, \mu)}+\frac{h(\must)}{l(\Xcst, \must)}.
\]
\end{theorem}

\proof
If $\Gammast$ has finite exponential growth relative to $(\Xcst, \dst)$ and $\must$ is non-elementary and of finite first moment,
then $(\partial \Xcst, \nu_\must)$ is a Poisson boundary for $(\Gammast, \must)$ \cite[Theorem 7.4]{Kaimanovich-hyperbolic}.
In particular, $h(\must)=h(\partial \Xcst, \must)$.
Therefore by Propositions \ref{Prop:upper_product} and \ref{Prop:lower_product},
if $l(\Xc, \mu)\ge l(\Xcst, \must)$,
then
for $\nu_\pi$-almost every $\xib \in \partial \Xc\times \partial \Xcst$,
\[
\lim_{r \to 0}\frac{\log \nu_\pi(\Bb(\xib, r))}{\log r}=\frac{h(\pi)-h(\must)}{l(\Xc, \mu)}+\frac{h(\must)}{l(\Xcst, \must)}.
\]
This shows that $\nu_\pi$ is exact dimensional.
The second claim follows from Lemma \ref{Lem:Frostman}.
\qed

\proof[Proof of Theorem \ref{Thm:product_intro}]
The claim follows from Theorem \ref{Thm:product} as a special case.
\qed

\begin{remark}\label{Rem:two_or_more}
Let us mention possible extensions and related questions.
\begin{enumerate}
\item The proof of Theorem \ref{Thm:product} can be extended to a product of more than two hyperbolic metric spaces.
For a positive $N \in \Z_+$,
let $(\Xc^{(i)}, d^{(i)})$ for $i=1, \dots, N$ be proper roughly geodesic hyperbolic metric spaces with bounded growth at some scale.
Further $\Gamma^{(i)}$ are countable subgroups in $\Isom \Xc^{(i)}$ with finite exponential growth for each $i=1, \dots, N$.
Let $\Gammab:=\Gamma^{(1)}\times \cdots \times \Gamma^{(N)}$.
For a probability measure $\pi$ on $\Gammab$ with non-elementary marginals $\mu^{(i)}$ of $\pi$ in $\Isom \Xc^{(i)}$ of finite first moment,
the harmonic measure $\nu_\pi$ on $\partial \Xc^{(1)} \times \cdots \times \partial \Xc^{(N)}$ is exact dimensional:
Let $\pi^{(i)}$ be the pushforward of $\pi$ to $\Isom \Xc^{(i)} \times \cdots\times \Isom \Xc^{(N)}$ for $i=1, \dots, N$.
If $l(\Xc^{(1)}, \mu^{(1)}) \ge \cdots \ge l(\Xc^{(N)}, \mu^{(N)})$,
then for $\nu_\pi$-almost every $\xib \in \partial \Xc^{(1)} \times \cdots \times \partial \Xc^{(N)}$,
\[
\dim \nu_\pi=
\sum_{i=1}^{N-1} \frac{h(\pi^{(i)})-h(\pi^{(i+1)})}{l(\Xc^{(i)}, \mu^{(i)})}+\frac{h(\pi^{(N)})}{l(\Xc^{(N)}, \mu^{(N)})}.
\]
In the above, $h(\pi^{(i)})$ denotes the asymptotic entropy for a $\pi^{(i)}$-random walk and $l(\Xc^{(i)}, \mu^{(i)})$ denotes the drift associated with a $\mu^{(i)}$-random walk for each $i=1, \dots, N$.
Further the Hausdorff dimension is computed by the quasi-metric defined as maximum of the ones in $\partial \Xc^{(i)}$.
The proof proceeds by the reverse induction in $i$ from $N$ to $1$ upon extending Propositions \ref{Prop:upper_product} and \ref{Prop:lower_product} and Theorem \ref{Thm:product} to the spaces $\Xc^{(i)} \times \cdots \times \Xc^{(N)}$ for $i=1, \dots, N$.
Since writing out all the details in this generality would hurt readability,
we refrain from producing the whole argument.
\item In \cite{Tdim}, the exact dimensionality of the harmonic measures for a single hyperbolic metric space $\Xc$ has been extended to several directions. 
For example, $\Xc$ can be replaced by a proper hyperbolic $\Xc$ without assuming bounded growth at some scale, and by a non-proper, separable and geodesic hyperbolic $\Xc$ with acylindrical action of a group.
In those cases, probability measures $\mu$ for $\mu$-random walks are assumed to satisfy that the support generates a non-elementary subgroup of isometries as a {\it semigroup} rather than a {\it group}.
It is expected that results in the present paper are extended to products of such hyperbolic spaces (with right assumption on random walks).
However, we need that the boundary be Polish (at least the space endowed with Borel structure be a standard Borel space) so that the conditional measures are well-defined.
Thus it is not clear as to whether the separability of $\Xc$ could be dropped.
\item It is not clear as to whether one can remove the condition on finite exponential growth relative to each factor in Theorem \ref{Thm:product}.
It is expected that the harmonic measure is exact dimensional without the assumption in regards of results in \cite{HochmanSolomyak} and \cite{LedrappierLessaExact}.
The issue in the present setting lies in a lack of Lebesgue differentiation theorem on boundaries.
This is available, for example, if the boundaries are Euclidean spaces, more generally, Riemannian manifolds, or if the harmonic measures are doubling (which is stringent).
In this paper, we have used a weaker version of Lebesgue differentiation theorem (Lemma \ref{Lem:BonkSchramm}).
However, we do not know if this would suffice to remove the condition on growth. 
See a related question in \cite[Quesion 4.3]{Tdim}.
\end{enumerate}
\end{remark}

\section{A positive lower bound for dimension}\label{Sec:positive}

\subsection{Pivotal times}\label{Sec:pivotal}
Let us recall the terminology and methods from \cite{Gouezel_Exp}.
For $\delta \in \R_+$, let $(\Xc, d)$ be a $\delta$-hyperbolic space with a base point $o$.
A sequence of points $x_0, \dots, x_n$ 
is called a {\bf $(C, D)$-chain} for some $C, D \in \R_+$
if 
\[
(x_{i-1}|x_{i+1})_{x_i} \le C \quad \text{for all $i=1, \dots, n-1$}, \quad
\text{and} \quad d(x_{i-1}, x_i)\ge D \quad \text{for all $i=1, \dots, n$}.
\]
If a sequence $x_0, \dots, x_n$ is a $(C, D)$-chain with $C \in \R_+$ and $D \ge 2C+2\delta+1$,
then 
\begin{equation}\label{Eq:Lem:3.7}
(x_0|x_n)_{x_1} \le C+\delta
\quad \text{and} \quad 
d(x_0, x_n) \ge \sum_{i=0}^{n-1}(d(x_i, x_{i+1})-(2C+2\delta))\ge n.
\end{equation}
See \cite[Lemma 3.7]{Gouezel_Exp}.
For $C \in \R_+$, $D=2C+2\delta+1$ and for $x, y \in \Xc$,
the {\bf chain shadow} $\Cc\Sc_x(y, C)$ of $y$ seen from $x$
is the set
\[
\Big\{z \in \Xc \ : \ \text{there exists a $(C, D)$-chain $x_0, \dots, x_n$; $x_0=x$, $x_n=z$ and $(x_0|x_1)_y \le C$} \Big\}.
\]

\begin{definition}\label{Def:Schottky}
For $\e, C, D \in \R_+$, 
a set of isometries $\Sc$ is called an $(\e, C, D)$-{\bf Schottky set} if the following three conditions are satisfied:
\begin{itemize}
\item[(1)] $\#\{s \in \Sc \ : \ (x|s\cdot y)_o \le C\} \ge (1-\e)\# \Sc$ for all $x, y \in \Xc$,
\item[(2)] $\#\{s \in \Sc \ : \ (x|s^{-1}\cdot y)_o \le C\} \ge (1-\e)\# \Sc$ for all $x, y \in \Xc$, and
\item[(3)] $d(o, s\cdot o) \ge D$ for all $s \in \Sc$.
\end{itemize}
\end{definition}

Let $\mu$ be a non-elementary probability measure on $\Isom \Xc$ with a countable support.
It is shown 
basically by a classical ping-pong argument
that for every $\e>0$
there exists a $C_0 \in \R_+$ satisfying the following:
For all $D \in \R_+$ there exists an $M \in \Z_+$ such that
the support of the $M$-fold convolution $\mu^{\ast M}$ of $\mu$ contains an $(\e, C_0, D)$-Schottky set $\Sc$ \cite[Corollary 3.13]{Gouezel_Exp}.
Let us fix the constants
\[
\e=1/100, \quad C_0 \in \R_+ \quad \text{and} \quad D \ge 20C_0+100\delta+1.
\]
Let $\lambda_\Sc$ denote the uniform distribution on $\Sc$.

Given a sequence of isometries $u_0, u_1, \dots$ on $\Xc$ and a sequence of independent random isometries $s_1, s_2, \dots$ with the identical distribution $\lambda_\Sc^{\ast 2}$,
let
\[
y_{n}^-:=u_0 s_1 u_1 \cdots s_{n-1} u_{n-1} \cdot o.
\]
Letting $s_i=a_i b_i$ where $a_i$ and $b_i$ are independent and distributed as $\lambda_\Sc$,
we define
\[
y_n:=u_0 s_1 u_1 \cdots s_{n-1}u_{n-1} a_n \cdot o
\quad \text{and} \quad y_n^+:=u_0 s_1 u_1 \cdots s_{n-1}u_{n-1} a_n b_n \cdot o.
\]
A sequence of {\bf pivotal times} $P_n \subseteq \{1, \dots, n\}$ is defined inductively as in the following:
Let $P_0:=\emptyset$ (the empty set).
Given $P_{n-1}$,
let $k:=\max P_{n-1}$ if $P_{n-1}\neq \emptyset$, and let $k=0$ and $y_k:=o$ if $P_{n-1}=\emptyset$.
Suppose that $P_{n-1}\neq \emptyset$.
Let us say that the {\bf local geodesic condition} is satisfied at time $n$ if 
\begin{equation}\label{Eq:lgc}
(y_k|y_n)_{y_n^-}\le C_0, \quad (y_n^-|y_n^+)_{y_n} \le C_0 \quad \text{and} \quad (y_n|y_{n+1}^-)_{y_n^+}\le C_0.
\end{equation}
If the local geodesic condition is satisfied at time $n$,
then we define
\[
P_n:=P_{n-1}\cup \{n\}.
\]
If otherwise, then letting $m$ be the largest pivotal time in $P_{n-1}$ such that
\[
y_{n+1}^- \in \Cc\Sc_{y_m}(y_m^+, C_0+\delta),
\]
we define 
$P_n:=P_{n-1}\cap \{1, \dots, m\}$, and $P_n:=\emptyset$ in the case when there is no such $m$. 
Note that the set $P_n$ depends only on the sequence $s_1, \dots, s_n$ for fixed $u_0, \dots, u_n$.

\begin{lemma}\label{Lem:4.4}
If $P_n:=\{k_1, \dots, k_p\}$ where $k_1<\cdots<k_p$,
then the sequence 
\[
o, y_{k_1}, y_{k_2}^-, y_{k_2}, \dots, y_{k_p}^-, y_{k_p}, y_{n+1}^-,
\]
forms a $(2C_0+4\delta, D-2C_0-3\delta)$-chain.
Moreover, if $D \ge 6 C_0+13\delta+1$,
then for every $i=2, \dots, p$,
the sequence $y_{k_i}^-, y_{k_i}, y_{n+1}^-$ is a $(2C_0+5\delta, D-6C_0-13\delta)$-chain.
\end{lemma}

\proof
The first claim is \cite[Lemma 4.4]{Gouezel_Exp}.
The second claim follows from the first claim and \eqref{Eq:Lem:3.7}.
Indeed, applying to them the sequence $y_{k_i}^-, y_{k_i}, \dots, y_{k_p}, y_{n+1}^-$
for each $i=2, \dots, p$ shows that
$(y_{k_i}^-|y_{n+1}^-)_{y_{k_i}} \le 2 C_0+5 \delta$,
and further, 
\[
d(y_{k_i}, y_{n+1}^-) \ge d(y_{k_i}, y_{k_{i+1}}^-)-2(2C_0+4\delta)-2\delta \ge D-6C_0-13\delta.
\]
The claim follows.
\qed

Let $\wbar s=(s_1, \dots, s_n)$.
Let us say that a sequence $\wbar s'=(s_1', \dots, s_n')$
where $s_i'=a_i'b_i'$ is {\bf pivoted from} $\wbar s$ if $\wbar s'$ and $\wbar s$ have the same pivotal times,
$b_i'=b_i$ for all $i=1, \dots, n$,
and $a_i'=a_i$ for all $i$ which is not a pivotal time.
The relation that $\wbar s'$ is pivoted from $\wbar s$ defines an equivalence relation among sequences.
In this notation, we understand that $b_i$ for all $i=1, \dots, n$ and $a_i$ for $i$ which is not a pivotal time are determined in $\wbar s$.

Let $\Ec_n(\wbar s)$ be the set of sequences which are pivoted from $\wbar s$.
Note that if $u_0, u_1, \dots$ are fixed,
then conditioned on $\Ec_n(\wbar s)$,
all $a_i$ are independent.
However, their distributions may depend on $i$.
For each $i=1, \dots, n$, let
\[
A_i(\wbar s):=\big\{a \in \Sc \ : \ \text{$s_i'=ab_i$ for some $\wbar s'=(s_1', \dots, s_i', \dots, s_n') \in \Ec_n(\wbar s)$} \big\}.
\]
It holds that for each pivotal time $i$ of $\wbar s$,
\[
\Pb(a_i=a \mid \Ec_n(\wbar s))=\lambda_\Sc(a\mid A_i(\wbar s))=\frac{\lambda_\Sc(a)}{\lambda_\Sc(A_i(\wbar s))} \quad \text{for $a \in A_i(\wbar s)$}.
\]
If $i$ is a pivotal time of $\wbar s$ and $\wbar s'=(s_1, \dots, s_i', \dots, s_n)$ in which $s_i=a_ib_i$ is replaced by $s_i'=a_i'b_i$ satisfies the local geodesic condition at time $i$,
then $\wbar s'$ is pivoted from $\wbar s$ \cite[Lemma 4.7]{Gouezel_Exp}.
By the definition of $(\e, C_0, D)$-Schottky set (Definition \ref{Def:Schottky}),
there are at most $2\e\#\Sc$ elements for which the local geodesic condition \eqref{Eq:lgc} does not hold at $i$.
Therefore for each pivotal time $i$ in $\wbar s$,
\begin{equation}\label{Eq:A_i}
\# A_i(\wbar s) \ge (1-2\e)\#\Sc.
\end{equation}

\begin{lemma}\label{Lem:inj}
Let $D \ge 10 C_0+25 \delta +1$.
For $\wbar s' \in \Ec_n(\wbar s)$, if $y_{n+1}^-=y_{n+1}'^-$, then $\wbar s=\wbar s'$.
\end{lemma}

\proof
Let $P_n=\{k_1, \dots, k_p\}$ with $k_1<\cdots<k_p$ be the set of pivotal times in $\wbar s$.
By Lemma \ref{Lem:4.4},
the sequences $y_{k_i}^-, y_{k_i}, y_{n+1}^-$ and $y_{k_i}'^-, y_{k_i}', y_{n+1}'^-$ are $(2C_0+5\delta, D-6 C_0-13\delta)$-chains respectively.
For $\wbar s' \in \Ec_n(\wbar s)$ with $\wbar s'\neq \wbar s$,
let $i$ be the first $i$ for which $s_{k_i} \neq s_{k_i}'$.
For $s_{k_i}=a_{k_i}b_{k_i}$ and $s_{k_i}'=a_{k_i}'b_{k_i}'$,
it holds that $a_{k_i}\neq a_{k_i'}$ and these $a_{k_i}$ and $a_{k_i}'$ are in the Schottky set $\Sc$,
whence $(a_{k_i}\cdot o|a_{k_i}'\cdot o)_o \le C_0$.
This shows that the sequence
\[
y_{n+1}', y_{k_i}', y_{k_i}^-, y_{k_i}, y_{n+1}^-, \quad \text{where $y_{k_i}^-=y_{k_i}'^-$},
\]
forms a $(2C_0+5\delta, D-6C_0-13\delta)$-chain.
For such $D$, one has $d(y_{n+1}'^-, y_{n+1}^-)>0$ by \eqref{Eq:Lem:3.7}, and thus $y_{n+1}'^-\neq y_{n+1}^-$, as required.
\qed

\subsection{A lower bound for entropy}

\begin{theorem}\label{Thm:entropy_lower_bound}
Let $\Gamma$ and $\Gammast$ be countable subgroups in $\Isom \Xc$ and in $\Isom \Xcst$ respectively, and $\Gammab:=\Gamma\times \Gammast$.
Further let us consider a probability measure $\pi$ on $\Gammab$ of the following form:
\[
\pi=\alpha \lambda \times \lambdast+(1-\alpha)\pi_0
\]
for some $\alpha \in (0, 1]$, non-elementary probability measures $\lambda$ and $\lambdast$ on $\Gamma$ and on $\Gammast$ respectively, and a probability measure $\pi_0$ on $\Gammab$.
Then for the marginal $\must$ of $\pi$ on $\Gammast$,
it holds that $h(\pi)-h(\partial \Xcst, \must)>0$.
Moreover, if in addition $\Gammab$ has a finite exponential growth relative to $(\Xc, \db)$,
then
for $\nu_{\must}$-almost every $\eta \in \partial \Xcst$,
the Hausdorff dimension of the conditional measure $\nu_\pi^\eta$ is positive.
\end{theorem}

Fix constants $\e=1/100$, $C_0\ge 0$ and $D \ge 20 C_0+100\delta+1$, and a $(1/100, C_0, D)$-Schottky set $\Sc$ in $\Gamma$, contained in the support of $\lambda^{\ast M}$ for some $M \in \Z_+$.
For $N:=2M$,
let us write for some $\beta \in (0, 1]$ and for a probability measure $\lambda_0$ on $\Gamma$,
\[
\lambda^{\ast N}=\beta \lambda_\Sc^{\ast 2}+(1-\beta)\lambda_0.
\]
Let us also write for a probability measure $\wbar\pi_0$ on $\Gammab=\Gamma\times \Gammast$,
\[
\pi^{\ast N}=\alpha^N \beta \lambda_\Sc^{\ast 2}\times \lambdast^{\ast N}+(1-\alpha^N \beta)\wbar \pi_0.
\]
For a sequence $\e_1, \e_2, \dots$ of independent Bernoulli random variables with the common parameter $\alpha^N \beta$,
let us define a sequence of independent random group elements $\gammab_1, \gammab_2, \dots$ where $\gammab_i=(\gamma_i, \gammast_i) \in \Gammab$ is distributed as $\lambda_\Sc^{\ast 2}\times \lambdast^{\ast N}$ if $\e_i=1$ and as $\wbar \pi_0$ if $\e_i=0$. 
Note in particular that for $i=1, 2, \dots$, conditioned on the event $\{\e_i=1\}$, random group elements $\gamma_i$ and $\gammast_i$ are independent.

Further we realize a $\pi$-random walk $\wb_{nN}$ at time $nN$ by a product $\gammab_1\cdots \gammab_n$ through a coupling on an enlarged probability space.
Let us define a sequence of pivotal times for $z_n=w_n\cdot o$ on $\Xc$,
where $\{w_n\}_{n \in \Z_+}$ is a $\mu$-random walk and $\mu$ is the marginal of $\pi$ on $\Gamma$.

Let $t_1, t_2, \dots$ be the sequence of $i$ such that $\e_i=1$.
For every positive $n \in \Z_+$,
let $\tau=\tau(n)$ be the maximum of $j$ with $Nt_j \le n$.
It holds that
\[
(N (t_j-1), Nt_j] \subset (0, n] \quad \text{for all $j=1, \dots, \tau$}.
\]
For each $j=1, \dots, \tau$,
let $s_{t_j}$ be $\gamma_j$, which is realized as the product of elements $x_i$ over $i\in (N (t_j-1), Nt_j]$ in the natural order from $\Z_+$.
Let us write $s_j':=s_{t_j}$ for brevity, and 
\[
u_0:=x_1\cdots x_{N (t_1-1)}, \quad u_j:=x_{Nt_j+1}\cdots x_{N (t_{j+1}-1)} \quad \text{and} \quad u(n):=x_{Nt_{\tau}+1}\cdots x_{n}.
\]
In the above, $u_0$, $u_j$ and $u(n)$ are defined as the identity if they are empty words.
For a $\mu$-random walk $w_n$ at time $n$, the orbit $z_n$ is realized as
\[
z_n=w_n\cdot o=u_0 s_1' u_1 \cdots s_\tau' u(n)\cdot o.
\]
Let $P_1, \dots, P_{\tau(n)}$ be the sequence of pivotal times of $z_n$ given $u_0, u_1, \dots, u(n)$.
Note that $P_{\tau(n)}$ depends not only on $\tau(n)$ but also on $n$.
It is shown that there exists a constant $\kappa>0$ such that
\begin{equation}\label{Eq:exp}
\Pb(\#P_{\tau(n)} \le \kappa n) \le e^{-\kappa n} \quad \text{for all $n \in \Z_+$},
\end{equation}
\cite[Proposition 4.11]{Gouezel_Exp}.

Note that the sequence $\{\e_i\}_{i=1}^\infty$ determines $\{t_i\}_{i=1}^\infty$ and $\tau=\tau(n)$ for each $n$.
Let us define the $\sigma$-algebra
\[
\Gc:=\sigma\Big(\text{$\e_i, t_i, \xst_i$ for $i=1,2, \dots$ and $x_i$ for $i \notin \bigcup_{j=1}^{\tau(n)}(N (t_j-1), Nt_j]$} \Big).
\]
Conditioning on $\Gc$ amounts to fixing a typical sequence $\{\e_i\}_{i=1}^\infty$, $\{t_i\}_{i=1}^\infty$ and a trajectory of $\must$-random walk $\{\wst_n\}_{n\in \Z_+}$, increments $x_i$ of $\mu$-random walk outside the time intervals $(N (t_j-1), Nt_j]$ for $j=1, \dots, \tau(n)$.
Under this conditioning, $s'_1, \dots, s'_\tau$ is a sequence of independent random elements in $\Gamma$ with the common distribution $\lambda_\Sc^{\ast 2}$.

Let $\wbar s:=(s'_1, \dots, s'_\tau)$, and $\Ec_\tau(\wbar s)$ be the set of sequences pivoted from $\wbar s$.
Conditioned on $\Ec_\tau(\wbar s)$ and $\Gc$,
the random group elements $a_i$ at pivotal times where $s'_i=a_i b_i$ are independent and each $a_i$ is distributed as $\lambda_\Sc(\,\cdot\mid A_i(\wbar s))$.
Let $\sigma(\Ec_\tau(\wbar s), \Gc)$ denote the $\sigma$-algebra generated by $\Ec_\tau(\wbar s)$ and $\Gc$.
Given $\Ec_\tau(\wbar s)$ and $u_0, \dots, u(n)$,
let us define the map
\[
\prod_{i \in P_{\tau(n)}} A_i(\wbar s) \to \Xc, \quad (a_i)_{i \in P_{\tau(n)}} \mapsto w_n\cdot o=u_0 a_1 b_1u_1\cdots u_{\tau-1}a_\tau b_\tau u(n)\cdot o.
\]
In the above, we understand that $\{b_i\}_{i=1, \dots, \tau(n)}$ and $\{a_i\}_{i \notin P_{\tau(n)}}$ are determined by $\Ec_\tau(\wbar s)$.
Under the conditioning, the map is injective by Lemma \ref{Lem:inj}.
Furthermore the conditional distribution $\Pb(w_n\cdot o \in \cdot\mid \sigma(\Ec_\tau(\wbar s), \Gc))$ is the pushforward by the injective map of the product measure $\lambda_\Sc(\,\cdot\mid A_i(\wbar s))$ over $i \in P_{\tau(n)}$ almost everywhere in $\Pb$.

\proof[Proof of Theorem \ref{Thm:entropy_lower_bound}]
Let us denote by $H(w)$ the entropy $H(\mu)$ for a random variable $w$ with the distribution $\mu$.
If $\Fc$ is a sub $\sigma$-algebra of $\sigma(\{(w_n, \wst_n)\}_{n \in \Z_+})$,
then the conditional entropy of $w_n$ with respect to $\Fc$ is defined by
\[
H(w_n\mid\Fc):=\Eb\Big[-\sum_{x \in \supp \mu_n}\Pb(w_n=x\mid \Fc)\log \Pb(w_n=x\mid \Fc)\Big].
\]
It holds that $H(w_n) \ge H(w_n \mid \Fc)$, further that if $\Fc_1$ and $\Fc_2$ are sub $\sigma$-algebras of $\sigma(\{(w_n, \wst_n)\}_{n \in \Z_+})$ and $\Fc_1 \subseteq \Fc_2$,
then $H(w_n \mid \Fc_1) \ge H(w_n \mid \Fc_2)$.
Since the $\sigma$-algebra $\sigma(\wst_n)$ generated by $\wst_n$ is included in $\Gc$, 
it follows that
\[
H(w_n\mid \sigma(\wst_n)) \ge H(w_n \mid \sigma(\Ec_\tau(\wbar s), \Gc)).
\]
Let us find a lower bound on the right hand side of the following:
\[
H(w_n\mid \sigma(\Ec_\tau(\wbar s), \Gc))=\Eb\Big[-\sum_{x \in \supp \mu_n}\Pb\(w_n=x\mid \sigma(\Ec_\tau(\wbar s), \Gc)\)\log \Pb\(w_n=x\mid \sigma(\Ec_\tau(\wbar s), \Gc)\)\Big].
\]
Since $\Pb(w_n\cdot o\in \cdot \mid \sigma(\Ec_\tau(\wbar s), \Gc))$ is the pushforward by an injective map of the product measure of $\lambda_\Sc(\cdot\mid A_i(\wbar s))$ over $i \in P_{\tau(n)}$,
one has $\Pb$-almost everywhere,
\[
-\sum_{x \in \supp \mu_n}\Pb(w_n=x\mid \sigma(\Ec_\tau(\wbar s), \Gc))\log \Pb(w_n=x\mid \sigma(\Ec_\tau(\wbar s), \Gc))=-\sum_{i \in P_{\tau(n)}}\log \frac{1}{\# A_i(\wbar s)}.
\]
This shows that for each $n \in \Z_+$ and for $\kappa>0$ in \eqref{Eq:exp},
\[
H(w_n \mid \sigma(\wst_n)) \ge \Eb\Big[\Big(\sum_{i \in P_{\tau(n)}}\log \# A_i(\wbar s)\Big)\cdot \1_{\{\#P_{\tau(n)}\ge \kappa n\}}\Big].
\]
If $i \in P_{\tau(n)}$,
then $\# A_i(\wbar s) \ge (1-2\e)\#\Sc$ by \eqref{Eq:A_i}
and $\Pb(\#P_{\tau(n)} \le \kappa n) \le e^{-\kappa n}$ by \eqref{Eq:exp},
the right hand side of the above inequality is at least
\[
\kappa n\log((1-2\e)\#\Sc)\cdot \Pb(\#P_{\tau(n)} \ge \kappa n) \ge \kappa n\log((1-2\e)\#\Sc)\cdot (1-e^{-\kappa n}).
\]
Therefore it holds that
\[
\liminf_{n \to \infty}\frac{1}{n}H(w_n\mid \sigma(\wst_n)) \ge \kappa \log((1-2\e)\#\Sc)>0.
\]
Noting that $H(w_n, \wst_n)=H(\wst_n)+H(w_n\mid\sigma(\wst_n))$ and $h(\partial \Xcst, \must) \le h(\must)$,
we obtain
\[
h(\pi)-h(\partial \Xcst, \must)\ge h(\pi)-h(\must)\ge \liminf_{n \to \infty}\frac{1}{n}H(w_n\mid \sigma(\wst_n))>0.
\]
(In the above, the second inequality is in fact the equality and the liminf is the limit.)
Thus the first claim follows.
The second claim follows from the first claim together with Theorem \ref{Thm:exact}.
\qed

\proof[Proof of Theorem \ref{Thm:entropy_lower_bound_intro}]
This is contained in Theorem \ref{Thm:entropy_lower_bound}.
\qed

\subsection*{Acknowledgments}
The author would like to thank the anonymous referee for a meticulous review and suggestions.
The author is partially supported by 
JSPS Grant-in-Aid for Scientific Research JP20K03602 and JP24K06711.

\bibliographystyle{alpha}
\bibliography{dc}

\end{document}